\documentclass[psamsfonts]{amsart}

\usepackage{graphicx}
\usepackage{amsmath}
\usepackage{amsfonts}
\usepackage{amssymb}
\usepackage{amsthm}
\usepackage{pstricks}
\usepackage{subfigure}
\usepackage{nicefrac}
\usepackage{tikz-cd}
\usetikzlibrary{decorations.pathmorphing}
\usepackage{url}
\usetikzlibrary{arrows.meta,decorations.pathreplacing}
\usetikzlibrary{arrows.meta}
\usepackage{float}
\input xy
\xyoption{all}

\newtheorem{theorem}{Theorem}[section]
\newtheorem{corollary}[theorem]{Corollary}
\newtheorem{lemma}[theorem]{Lemma}
\newtheorem{proposition}[theorem]{Proposition}
\newtheorem{question}[theorem]{Question}

\theoremstyle{definition}
\newtheorem{definition}[theorem]{Definition}
\newtheorem{example}[theorem]{Example}

\newtheorem{remark}[theorem]{Remark}

\title[Generalized Hopf bifurcations: homology and obstructions]{Generalized Hopf bifurcations: a homological analysis and obstructions to their existence}

\author{H\'ector Barge}
\address{E.T.S. Ingenieros Inform\'aticos, Universidad Polit\'ecnica de Madrid, 28660 Madrid, Espa\~na}
\curraddr{}
\email{h.barge@upm.es}
\thanks{The authors are is partially supported by the Spanish Ministerio de Ciencia e Innovaci\'on grant PID2024-156663NB-I00 and the first author is also supported by grant PID2025-174633NA-I00}

\author{J.J. S\'anchez-Gabites}
\address{Facultad de Ciencias Matem\'aticas, Universidad Complutense de Madrid, 28040 Madrid, Espa\~na}
\curraddr{}
\email{jajsanch@ucm.es}

\subjclass[2020]{Primary: 37G35, 37B35; Secondary: 57N99}
\keywords{Attractor, Generalized Hopf Bifurcation, Basin of Attraction, Topological Spine}
\date{}

\begin{document}

\maketitle

\begin{abstract}
    A generalized Hopf bifurcation takes place when varying a parameter of a (semi)dynamical system causes an attractor $K$ to become a repeller. In doing so, it expels a (parameter-dependent) family of attractors.  In this paper we concentrate on two questions: (1) what do the expelled attractors look like?, and (2) can any attractor $K$ undergo a generalized Hopf bifurcation? Regarding (1), we compute the cohomology of the expelled attractors in terms of the cohomology of $K$; in the case of continuous-time dynamics this can be improved to a computation of their Borsuk shape. As for (2), we show that certain compacta cannot undergo generalized Hopf bifurcations and, when the ambient space has dimension $ \leq 3$, give a complete topological characterization of which can.
\end{abstract}

\section{Introduction}

A well-known phenomenon in smooth planar dynamics is when an attracting fixed point becomes, when varying a parameter, a repelling fixed point. Under appropriate analytical assumptions, the fixed point expels a closed invariant curve as a consequence of this change of stability. This kind of bifurcation was originally studied by Poincar\'e, Andronov and Hopf in the case of flows and by Naimark and Sacker in the case of diffeomorphisms. It has subsequently been extended to bifurcations of periodic orbits and also in higher dimensional phase spaces. For the sake of brevity we will refer to all these as classical Hopf bifurcations.

From the point of view of Topological Dynamics there is a natural generalization of this setup as follows. The phase space is a boundariless manifold $M$ where we have a semidynamical system that depends on a parameter $\lambda \in [0,1]$. Thus, the dynamics is generated by a family of embeddings $f_{\lambda} : M \longrightarrow M$ or injective semiflows $\varphi_{\lambda}$ on $M$. Following Seibert and Florio (\cite{MR1277591}), we say that a compact set $K \subseteq M$ undergoes a generalized Hopf bifurcation if it is invariant for all $\lambda$, it is an attractor for $\lambda = 0$ and becomes a repeller for $\lambda > 0$. The reason to consider semidynamical systems is motivated in detail in Section \ref{sec:hopf0}; here we just mention that it allows to accommodate the frequent situation where the change in stability of $K$ takes place only inside a submanifold of ambient space. Generalized Hopf bifurcations extend the classical setup in two directions: (i) the set $K$ can be completely arbitrary, rather than a critical point or a periodic orbit, and (ii) there are no analytic assumptions on how the transition from $K$ being an attractor to a repeller happens. It can easily be shown that there appears a family of attractors $A_{\lambda}$, all disjoint from $K$, that converge upper semicontinuously to $K$ as $\lambda \rightarrow 0^+$. These are called the attractors ``expelled'' by $K$. For instance, in the classical Hopf bifurcation in the plane, these $A_{\lambda}$ would be simple closed curves (periodic orbits in the case of a flow). In the generalized setting there is no reason to expect that the $A_{\lambda}$ have a nice topological structure, and this is suggested by Figure~\ref{fig:hopf_bifurcation_compact_clean}.

Two natural questions arise:
\begin{enumerate}
     \item What do the expelled attractors look like?
    \item What compacta $K$ admit this kind of bifurcations?
\end{enumerate}

\begin{figure}[htbp]
\centering
\begin{tikzpicture}[scale=0.95, >=Stealth]

  \def\Kpath{(0.47, 0.00) -- (0.79, 0.14) -- (0.55, 0.18) -- (0.47, 0.28) -- (0.34, 0.39) -- (0.14, 0.39) -- (0.02, 0.50) -- (-0.22, 0.45) -- (-0.29, 0.31) -- (-0.37, 0.25) -- (-0.70, 0.25) -- (-0.52, 0.09) -- (-0.77, 0.02) -- (-0.72, -0.14) -- (-0.52, -0.21) -- (-0.57, -0.35) -- (-0.23, -0.23) -- (-0.11, -0.27) -- (0.02, -0.44) -- (0.18, -0.48) -- (0.23, -0.25) -- (0.36, -0.23) -- (0.57, -0.23) -- (0.44, -0.09) -- cycle}

  \begin{scope}[xshift=0cm, yshift=0cm]


    \filldraw[gray!50, thick] \Kpath;
    \node[font=\small] at (0,0) {$K$};

    \draw[->, semithick, gray!80] (1.8, 1.0) to[out=-90, in=20] (0.65, 0.25);
    \draw[->, semithick, gray!80] (-1.8, 1.0) to[out=-90, in=160] (-0.65, 0.25);
    \draw[->, semithick, gray!80] (1.8, -1.0) to[out=90, in=-20] (0.60, -0.25);
    \draw[->, semithick, gray!80] (-1.8, -1.0) to[out=90, in=-160] (-0.60, -0.25);

    \node[font=\small, align=center] at (0, -2.5) {(a) $\lambda = 0$};

  \end{scope}

  \begin{scope}[xshift=7.0cm, yshift=0cm]

    \draw[thick, blue!70!black, decorate, decoration={random steps, segment length=4pt, amplitude=3pt}] (0,0) ellipse (2.0cm and 1.25cm);

    \filldraw[gray!50, thick] \Kpath;
    \node[font=\small] at (0,0) {$K$};

    \node[font=\small, text=blue!80!black] at (0, 1.6) {$A_\lambda$};

    \draw[->, semithick, red!75!black] (0.65, 0.25) to[out=25, in=-135] (1.48, 0.70);
    \draw[->, semithick, red!75!black] (-0.65, 0.25) to[out=155, in=-45] (-1.48, 0.70);
    \draw[->, semithick, red!75!black] (0.60, -0.25) to[out=-25, in=135] (1.48, -0.70);
    \draw[->, semithick, red!75!black] (-0.60, -0.25) to[out=-155, in=45] (-1.48, -0.70);

    \draw[->, semithick, black!70] (2.5, 1.5) to[out=-110, in=35] (1.72, 0.90);
    \draw[->, semithick, black!70] (-2.5, 1.5) to[out=-70, in=145] (-1.72, 0.90);
    \draw[->, semithick, black!70] (2.5, -1.5) to[out=110, in=-35] (1.72, -0.90);
    \draw[->, semithick, black!70] (-2.5, -1.5) to[out=70, in=-145] (-1.72, -0.90);

    \node[font=\small, align=center] at (0, -2.5) {(b) $\lambda > 0$};

  \end{scope}

\end{tikzpicture}
\caption{A generalized Hopf bifurcation of a point $K$. For $\lambda = 0$ it is an attractor; for $\lambda > 0$ it becomes a repeller and expels an attractor $A_{\lambda}$}
\label{fig:hopf_bifurcation_compact_clean}
\end{figure}

There is already a fair amount of literature dealing with question (1). For instance, in \cite{MR2371787} the author describes, in the context of flows, the Borsuk shape of the attractors expelled by the bifurcation of a fixed point or a periodic trajectory. There are also some results for discrete dynamical systems, but they are only for fixed points in the plane \cite{Barge2020BifurcationsSystems} or, in higher dimensions, impose quite strong assumptions on $K$ and the phase space \cite{MR4421505}. Because of its relevance in population dynamics, the bifurcation of fixed points in the boundary of a manifold has also been studied in \cite{MR3199752}. In this paper we give a complete homological description of the expelled attractors for any subcompactum of a boundariless manifold. This is done in the most general context of discrete semidynamical systems, and so it encompasses most of the existing literature. Since the expelled attractors may have a very complicated topology, the preferred homological language to capture their global features is \v{C}ech cohomology, rather than singular homology, and this is what we will use throughout the paper. In the case of continuous-time dynamics (i.e. injective semiflows) we give an even more detailed description of the expelled attractors by using Borsuk's theory shape, which is a sort of homotopy theory especially suited to complicated compacta.

As for question (2), as far as we know it has not been studied before in the literature. Its precise meaning is this: given a compactum $K$ in a boundariless manifold $M$, does there exist a (parameterized) semidynamical system for which $K$ undergoes a generalized Hopf bifurcation? Obviously, a necessary condition for this to happen is that $K$ can be realized as an attractor in the first place. In the context of continuous-time dynamics, we show that this is also a sufficient condition. By contrast, this is no longer true in the discrete-time case. Indeed, we will prove that a compactum that undergoes a generalized Hopf bifurcation must have a finitely generated \v{C}ech cohomology. Thus, for example, the standard dyadic solenoid is an attractor for a homeomorphism of $\mathbb{R}^3$ but does not admit a Hopf bifurcation because its cohomology is not finitely generated. When $M$ is an orientable $3$-manifold we obtain a full answer to question (2): a compactum $K$ in $M$ admits a Hopf bifurcation if and only if it is a topological spine. This latter condition is concerned with how $K$ sits in $M$ and, in a certain sense, is very restrictive among attractors of discrete dynamical systems. Heuristically, this has the consequence that the topology of generalized Hopf bifurcations in dimension $3$ is the same both for discrete and for continuous dynamics. This characterization also shows that the main result of \cite{MR4421505}, where $K$ is required to satisfy an assumption equivalent to being a topological spine, actually holds in full generality in orientable $3$-manifolds. 

The two questions stated above are much more difficult to study for embeddings (i.e. discrete-time dynamics) compared to semiflows (i.e. continuous-time dynamics). Roughly, the reason is that semiflows provide homotopy equivalences in a natural way whereas in the discrete-time setting one has to replace those with algebraic arguments. Thus, most of the paper is set up for embeddings, and only in the last section we turn to semiflows. The paper is organized as follows. In Section~\ref{sec:completion} we describe how to complete the basin of attraction of an attractor of an embedding by extending ``backwards'' both the phase space and the embedding itself, so that the latter becomes a homeomorphism. This construction removes ``artificial holes'' in the basin of attraction and allows us to relate the topology of the attractor and its completed basin of attraction. In Section~\ref{sec:2} we show that the completion of the basin of attraction is rigid, i.e., its topological type is invariant under continuation. This result, which may be of independent interest, is a key tool for the remaining of the paper. In Section~\ref{sec:hopf0} we discuss the motivation of the generalized Hopf bifurcation in more detail and show that attractors of injective semiflows always admit this kind of bifurcations. In Section~\ref{sec:homological} we prove that if a compactum admits a generalized Hopf bifurcation it must have a finitely generated \v Cech cohomology. We also determine the \v Cech cohomology of the expelled attractors in terms of that of $K$. In Section~\ref{prop:spinesattractors} we relate the dynamical property of admitting a generalized Hopf bifurcation with the topological property of being a topological spine. We see that these are equivalent in the continuous-time case and also in the discrete-time case when $M$ is an orientable $3$-manifold. Finally, in Section~\ref{sec:semi} we show that the completion of the basin of attraction can be also carried out for injective semiflows and used it to determine the Borsuk's shape of the expelled attractors.

\section{Attractors of embeddings and completing their basins of attraction} \label{sec:completion}

Let $M$ be a locally compact metric space, and $f:M\longrightarrow M$ an embedding (i.e. a continuous and injective map). We say that $K\subseteq M$ is an \emph{attractor} if it has a compact neighbourhood $P$ such that $f(P)\subseteq {\rm int}(P)$ and
\[
K=\bigcap_{n\geq 0} f^n(P).
\]
Such a neighbourhood $P$ is called a \emph{trapping region} for $K$. It is not difficult to see that attractors have the following properties (which in fact characterize them without any reference to trapping regions):
\begin{enumerate}
    \item $K$ is compact, $f(K)=K$ and $f|_{K}$ is a homeomorphism.
    \item $K$ attracts nearby points, that is, it has a neighbourhood $U$ such that for every $x$ in $U$, ${\rm dist}(f^n(x),K) \rightarrow 0$ as $n \rightarrow +\infty$
    \item $K$ is stable in the sense of Lyapunov, that is, every neighbourhood $U$ of $K$ contains a neighbourhood $V$ of $K$ such that $f^n(V)\subseteq U$ for each $n\geq 0$.
\end{enumerate}

Given an attractor $K$, the set $\mathcal{A}(K)$ of points attracted by $K$ is called its \emph{basin of attraction}. Since the forward iterates of any such point must eventually fall in ${\rm int}\ P$ we have that $\mathcal{A}(K) = \bigcup_{n \geq 0} f^{-n}({\rm int}\ P)$, which shows that $\mathcal{A}(K)$ is an open neighbourhood of $K$. Lyapunov stability ensures that not only points, but in fact compact subsets of $\mathcal{A}(K)$, are (uniformly) attracted by $K$. Explicitly, if $B \subseteq \mathcal{A}(K)$ is compact, for any neighbourhood $V$ of $K$ there exists $n_0$ such that $f^n(B) \subseteq V$ for all $n \geq n_0$. We will use this property tacitly in later sections.

Clearly $\mathcal{A}(K)$ is forward invariant, so we have a well defined restriction $f : \mathcal{A}(K) \longrightarrow \mathcal{A}(K)$. It is injective and, near $K$, also surjective. Points far enough from $K$ may not have a preimage, but there is a natural way of supplying these by abstractly pulling backwards the basin of attraction. We describe this now. 

Let $A_0, A_1, \ldots$ be countably many disjoint copies of $\mathcal{A}(K)$ and take their union $\coprod A_i$ topologized in the usual way. Identify each $x \in A_i$ with $f(x) \in A_{i+1}$ (and then in turn also with $f^2(x) \in A_{i+2}, \ldots$), obtaining a quotient space $\hat{\mathcal{A}}(K)$. The reader familiar with direct limits may recognize this as the direct limit of the direct sequence \[
\begin{tikzcd}
\mathcal{A}(K) \arrow[r, "f"] & \mathcal{A}(K)\arrow[r, "f"] &\ldots
\end{tikzcd}
\]

Given a point $x \in \mathcal{A}(K)$, we may regard it as belonging to any one of the $A_i$ and then consider its equivalence class in $\hat{\mathcal{A}}(K)$. We denote this class by $[x \in A_i]$. Thus by definition $[x \in A_i] = [y \in A_j]$ if and only if there is $k \geq i,j$ such that $f^{k-i}(x) = f^{k-j}(y)$. There is a well-defined continuous map $\hat{f} : \hat{\mathcal{A}}(K) \longrightarrow \hat{\mathcal{A}}(K)$ given by $[x \in A_i] \longmapsto [f(x) \in A_i]$. It has a continuous inverse given by $[x \in A_i] \longmapsto [x \in A_{i+1}]$, and so $\hat{f}$ is a homeomorphism of $\hat{\mathcal{A}}(K)$. The following proposition justifies that we regard $(\hat{\mathcal{A}}(K), \hat{f})$ as a completion of the basin of attraction of $K$ and its dynamics. (Perhaps the word ``extension'' would be more appropriate, but that already has a different meaning in Dynamics).

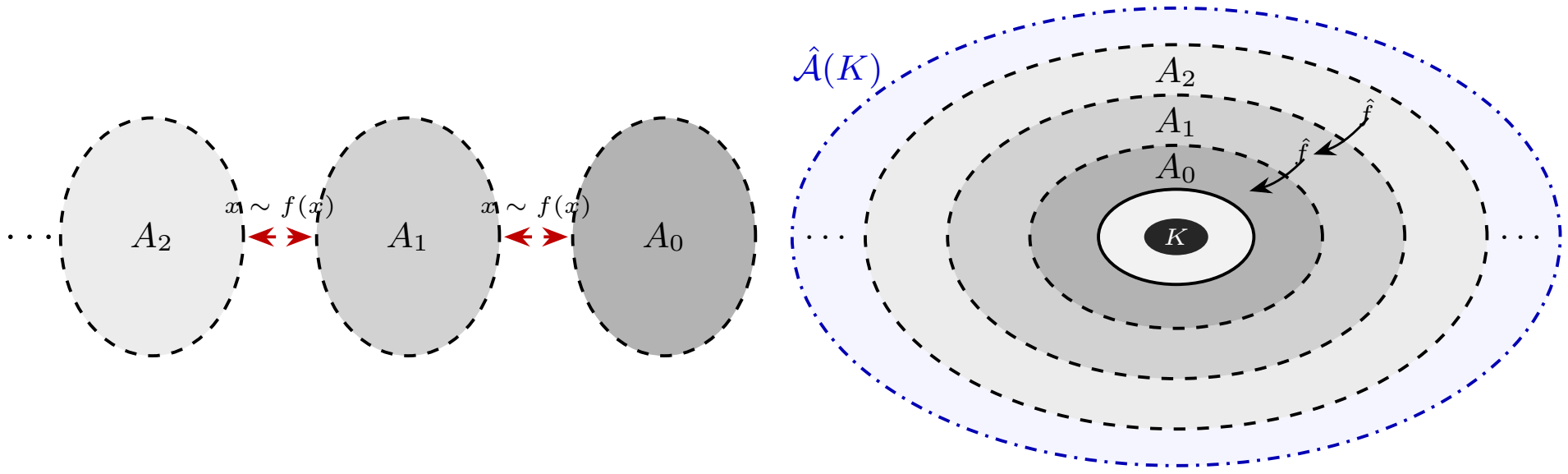
\begin{figure}[htbp]
\centering
\resizebox{\textwidth}{!}{
\begin{tikzpicture}[scale=0.82, >=Stealth]

  \node[font=\large] at (-0.3, 0) {$\dots$};

  \fill[gray!15] (1.0,0) ellipse (1.0cm and 1.3cm);
  \draw[thick, dashed] (1.0,0) ellipse (1.0cm and 1.3cm);
  \node[font=\small] at (1.0, 0) {$A_2$};

  \fill[gray!35] (3.8,0) ellipse (1.0cm and 1.3cm);
  \draw[thick, dashed] (3.8,0) ellipse (1.0cm and 1.3cm);
  \node[font=\small] at (3.8, 0) {$A_1$};

  \fill[gray!60] (6.6,0) ellipse (1.0cm and 1.3cm);
  \draw[thick, dashed] (6.6,0) ellipse (1.0cm and 1.3cm);
  \node[font=\small] at (6.6, 0) {$A_0$};

  \draw[<->, red!75!black, thick, dashed] (2.05, 0) -- node[above=1pt, font=\tiny, text=black] {$x \sim f(x)$} (2.75, 0);
  \draw[<->, red!75!black, thick, dashed] (4.85, 0) -- node[above=1pt, font=\tiny, text=black] {$x \sim f(x)$} (5.55, 0);

  \begin{scope}[xshift=12.2cm, yshift=0cm]

    \fill[blue!4] (0,0) ellipse (4.2cm and 2.5cm);
    \draw[thick, dash dot, blue!70!black] (0,0) ellipse (4.2cm and 2.5cm);
    \node[blue!80!black, font=\small\bfseries] at (-3.7, 1.85) {$\hat{\mathcal{A}}(K)$};

    \node[font=\small] at (3.8, 0) {$\dots$};
    \node[font=\small] at (-3.8, 0) {$\dots$};

    \fill[gray!15, even odd rule] (0,0) ellipse (3.4cm and 2.1cm) (0,0) ellipse (2.5cm and 1.55cm);
    \fill[gray!35, even odd rule] (0,0) ellipse (2.5cm and 1.55cm) (0,0) ellipse (1.6cm and 1.0cm);
    \fill[gray!60, even odd rule] (0,0) ellipse (1.6cm and 1.0cm) (0,0) ellipse (0.85cm and 0.52cm);

    \fill[gray!10] (0,0) ellipse (0.85cm and 0.52cm);
    \fill[black!85] (0,0) ellipse (0.35cm and 0.2cm);
    \node[white, font=\tiny] at (0,0) {$K$};

    \draw[thick, dashed] (0,0) ellipse (3.4cm and 2.1cm);
    \draw[thick, dashed]         (0,0) ellipse (2.5cm and 1.55cm);
    \draw[thick, dashed] (0,0) ellipse (1.6cm and 1.0cm);
    \draw[thick]         (0,0) ellipse (0.85cm and 0.52cm);

    \node[font=\small] at (0, 1.80) {$A_2$};
    \node[font=\small] at (0, 1.25) {$A_1$};
    \node[font=\small] at (0, 0.75) {$A_0$};

    \draw[->, semithick] (2.1, 1.3) to[bend left=15] node[above right, font=\tiny] {$\hat{f}$} (1.5, 0.9);
    \draw[->, semithick] (1.4, 0.85) to[bend left=15] node[above right, font=\tiny] {$\hat{f}$} (0.8, 0.5);

  \end{scope}

\end{tikzpicture}%
} 
\caption{Construction of the completion $\hat{\mathcal{A}}(K)$.}
\label{fig:abstract_basin_gluing}
\end{figure}

\begin{proposition} \label{prop:extendbasin} Suppose the embedding $f$ is an open map. Then there is a homeomorphic embedding $e$ of $\mathcal{A}(K)$ as an open invariant subset of $\hat{\mathcal{A}}(K)$. This embedding satisfies $\hat{f} \circ e = e \circ f$. The compact set $e(K)$ is an attractor for $\hat{f}$ whose basin of attraction is all of $\hat{\mathcal{A}}(K)$. 
\end{proposition}
\begin{proof} Let $e : \mathcal{A}(K) \longrightarrow \hat{\mathcal{A}}(K)$ be defined by $x \longmapsto [x \in A_0]$. The fact that $f$ is injective readily implies that $e$ is injective as well, and it is also straightforward to check that $\hat{f} \circ e = e \circ f$. We claim that $e$ is actually a homeomorphism onto an open subset of $\hat{\mathcal{A}}(K)$. Let $U$ be an open subset of $\mathcal{A}(K)$. To check that $e(U)$ is open in $\hat{\mathcal{A}}(K)$ we need to show that the set of points $x \in A_i$ such that $[x \in A_i] \in e(U)$ is open in $A_i$. But that set is readily seen to be $f^i(U)$. Since $f$ is an open map, the $f^i(U)$ are all indeed open. 

If $P$ is a trapping region for $K$, clearly $e(P)$ is a trapping region for $e(K)$ because \[\hat{f}(e(P)) = e(f(P)) \subseteq e({\rm int}_{\mathcal{A}(K)}\ P) = {\rm int}_{\hat{\mathcal{A}}(K)}\ e(P),\] where we have used that $e$ is an open map. Thus $e(K)$ is an attractor for $\hat{f}$ in $\mathcal{\hat{A}}(K)$. To check that every point in $\hat{\mathcal{A}}(K)$ is attracted by $K$ it suffices to show that it eventually enters $e(P)$. Let $[x \in A_i]$ and pick $n$ so large that $f^{n-i}(x) \in P$. In particular $f^n(x) \in f^i(P)$ so $f^n(x) = f^i(y)$ for some $y \in P$, and $\hat{f}^n([x \in A_i]) = [f^n(x) \in A_i] = [f^i(y) \in A_i] = [y \in A_0] \in e(P)$.
\end{proof}

\begin{remark} \label{rem:completion} (1) If $f$ was a homeomorphism to begin with, then $e$ is actually onto and so the completion $\mathcal{\hat{A}}(K)$ is just (a homeomorphic copy of) $\mathcal{A}(K)$ itself. In general, $\hat{\mathcal{A}}(K)$ is larger than $\mathcal{A}(K)$ and need not be (embeddable as) a subset of $M$.

(2)  From now on we remove $e$ from the notation and just regard $\mathcal{A}(K)$ as an open invariant subset of the completion $\hat{\mathcal{A}}(K)$. The same convention will be followed with $K$, small neighbourhoods of $K$, etc. since these all lie in $\mathcal{A}(K)$.

(3) We will only be interested in the case when $M$ is a boundariless manifold. Then $f$ is automatically an open map by the theorem on invariance of domain, so Proposition \ref{prop:extendbasin} applies. Also, $\hat{\mathcal{A}}(K)$ is a boundariless manifold as well. To check this notice first that $\mathcal{A}(K)$ is a boundariless manifold because it is open in $M$. Since $K$ is a global attractor in $\hat{\mathcal{A}}(K)$ and $\mathcal{A}(K)$ is an open neighbourhood of $K$, we have that $\hat{\mathcal{A}}(K)$ is the union of the increasing sequence of open sets $\hat{f}^{-n}(\mathcal{A}(K))$. Each of these is a boundariless manifold because it is homeomorphic to $\mathcal{A}(K)$ via $\hat{f}^{-n}$. Thus their union $\hat{\mathcal{A}}(K)$ is a boundariless manifold.

(4) Continuing with (3), suppose in addition that $M$ is orientable along $K$. This means that $K$ has a small open neighbourhood (hence a boundariless manifold) $U$ in $M$ that is a orientable. Then $\hat{\mathcal{A}}(K)$ is orientable. Otherwise $\hat{\mathcal{A}}(K)$ would contain an orientation-reversing loop $\gamma$, but then $\hat{f}^n(\gamma)$ would be an orientation-reversing loop for all $n$, and for $n$ large enough $\hat{f}^n(\gamma) \subseteq U$, contradicting the orientability of $U$.
\end{remark}

We record the following for later use:

\begin{remark} \label{rem:completionP} Here is a different description of $\hat{\mathcal{A}}(K)$. Let $P$ be a trapping region for $f$. We can repeat the same direct limit construction for $P$: we form the disjoint union $\coprod P_i$ of countably many copies $P_i$ of $P$ and quotient by the equivalence relation generated by $x \in P_i \sim f(x) \in P_{i+1}$. Denote by $\hat{\mathcal{A}}_P(K)$ the resulting quotient space. Clearly there is a continuous map $h : \hat{\mathcal{A}}_P(K) \longrightarrow \hat{\mathcal{A}}(K)$ defined by $h([x \in P_i]) := [x \in A_i]$. This map has an inverse $g$ defined as follows: given $[x \in A_i]$, find $n = n(x)$ so large that $f^{n} \in {\rm int}\ P$ and set $g([x \in A_i]) := [f^{n}(x) \in P_{i+n}]$. This is independent of the choice of $n$ since replacing $n$ by $n+1$ produces the same point in $\hat{\mathcal{A}}_P(K)$. Also, $g$ is continuous map: $x$ has a neighbourhood $U$ such that $f^{n}(U) \subseteq {\rm int}\ P$, so for all points $y$ in $U$ we may take the same $n(y) = n$ and therefore locally $g$ is just $f^n$. Therefore $(\hat{\mathcal{A}}_P(K),\hat{f}_P)$ is conjugated to $(\hat{\mathcal{A}}(K),\hat{f})$.
\end{remark}

Later on we will need to exploit a relation between the homology of an attractor and its basin of attraction or, more precisely, between the \v{C}ech cohomology $\check{H}^*(K)$ of an attractor $K$ and the (singular) cohomology $H^*(\hat{\mathcal{A}}(K))$ of its completed basin of attraction. We first state the result, which is a modification of Theorem~2 from \cite{MR3249272}, and then discuss it in some detail.

\begin{theorem}\label{thm:pacojaime}
Let $M$ be a boundariless manifold. For technical reasons, assume $M$ is smooth or triangulable. Let $K\subseteq M$ be an attractor of an embedding $f:M\longrightarrow M$ and $\hat{\mathcal{A}}(K)$ the completion of its basin of attraction. Then the following are equivalent:
\begin{itemize}
	\item[(i)] The inclusion $i:K\hookrightarrow\hat{\mathcal{A}}(K)$ induces isomorphisms $i^*: H^*(\hat{\mathcal{A}}(K))\longrightarrow\check{H}^*(K)$.
	\item[(ii)] The cohomology $\check{H}^*(K)$ is finitely generated.
	\item[(iii)] The cohomology $H^*(\hat{\mathcal{A}}(K))$ is finitely generated.
\end{itemize}
\end{theorem}

The reader may be unfamiliar with \v{C}ech cohomology. It is a cohomology theory which coincides with the usual cohomology over nice spaces (say, manifolds or polyhedra) but has several advantages over the latter when applied to compacta that have a complicated topology, as is usually the case with attractors. This is why the theorem above involves the \v{C}ech cohomology of $K$. The reason why it is related to the cohomology of $\hat{\mathcal{A}}(K)$ and not necessarily to that of $\mathcal{A}(K)$ can be easily illustrated with an example. Take $K = \{0\} \subseteq \mathbb{R}^2$ and the homeomorphism $f(x) = x/2$. Let $M$ be the open unit disk with a few disjoint closed disks $D_i$ removed. The $D_i$ are required to lie entirely in the annulus $1/2 < \|x\| < 1$. Then clearly $f|_M : M \longrightarrow M$ is an embedding and $K$ is an attractor whose basin of attraction $\mathcal{A}(K)$ is all of $M$. Evidently the cohomologies of $K$ and $\mathcal{A}(K)$ do not agree, even though the cohomology of $K$ is finitely generated. However, here $(\hat{\mathcal{A}}(K),\hat{f})$ is essentially $\mathbb{R}^2$ equipped with the original $f$, and obviously the cohomologies of $K$ and $\hat{\mathcal{A}}(K)$ do agree. This provides a different motivation for the introduction of the completion $\hat{\mathcal{A}}(K)$ of the basin of attraction: in some sense, it fills in any ``artificial holes'' that $\mathcal{A}(K)$ may have.

The proof of Theorem \ref{thm:pacojaime} is fairly technical, but we can say a few words here (full details are given in \cite{MR3249272}). Let $P$ be a compact neighbourhood of $K$ contained in its basin of attraction. By replacing $f$ with some iterate we may take $P$ to be positively invariant. Consider the map $f^* : H^*(P) \longrightarrow H^*(P)$ and assume for the sake of argument that its succesive images, which form a decreasing sequence of subgroups of $H^*(P)$, eventually stabilize to some subgroup $G$ (roughly, this is ensured by any one of the equivalent conditions of the theorem). One readily sees that this implies that the direct limit of \begin{equation} \label{eq:dirlim} \xymatrix{H^*(P) \ar[r]^{f^*} & H^*(P) \ar[r]^-{f^*} & \ldots } \end{equation} and the inverse limit of \begin{equation} \label{eq:invlim} \xymatrix{H^*(P) & \ar[l]_{f^*} H^*(P) & \ar[l]_-{f^*} \ldots }\end{equation} are both isomorphic to $G$. The attractor $K$ can be identified with the inverse limit of the inverse sequence $\xymatrix{P & \ar[l]_{f} P & \ar[l]_{f} \ldots}$ and a crucial property of \v{C}ech cohomology, called continuity, implies that $\check{H}^*(K)$ is given by the direct limit of \eqref{eq:dirlim}. On the other hand, by definition $\hat{\mathcal{A}}(K)$ is the direct limit of the sequence  $\xymatrix{P \ar[r]^f & P \ar[r]^f & \ldots}$ but, since cohomology does not always behave well with respect to direct limits, in principle $H^*(\hat{\mathcal{A}}(K))$ may not be the inverse limit of \eqref{eq:invlim}. Most of the work in \cite{MR3249272} is devoted to an analysis of this difficulty, eventually showing that indeed the cohomology of $\hat{\mathcal{A}}(K)$ is the inverse limit of \eqref{eq:invlim} and so $\check{H}(K) = G = H^*(\hat{\mathcal{A}}(K))$ as desired. That analysis makes use of a certain algebraic tool due to Milnor that needs $P$ to have good homotopy extension properties. Simple examples of such $P$ are manifolds or CW-complexes. This is where we use the technical assumption that the manifold $M$ be smooth or triangulable: in that case it is straightforward to show that any compact set $K \subseteq M$ has a neighbourhood basis of compact manifolds, for instance.

\begin{remark} \label{rem:pacojaime} If one takes cohomology with coefficients in a field $\mathbb{F}$, the three items in Theorem \ref{thm:pacojaime} hold automatically: the inclusion induced map $H^*(\hat{\mathcal{A}}(K);\mathbb{F}) \longrightarrow \check{H}^*(K;\mathbb{F})$ is an isomorphism and both cohomology groups are finitely generated; i.e. finite dimensional vector spaces.
\end{remark}

\section{Invariance of the topology of the basin of attraction under perturbations}\label{sec:2}

The setup for this section is as follows. The ambient space $M$ is a boundariless (topological) manifold. We have a family of embeddings $f_{\lambda}:M\longrightarrow M$ continuously parametrized by $\lambda\in[0,1]$. This means that the map $F:M\times[0,1]\longrightarrow M\times[0,1]$ given by
\[
F(x,\lambda)=(f_\lambda(x),\lambda)
\]
is continuous.

Assume that for a certain $\lambda_0\in [0,1]$ the embedding $f_{\lambda_0}$ has a trapping region $P$. By the continuous dependence on the parameter $\lambda$, there exists some neighbourhood $J$ of $\lambda_0$ in $[0,1]$ such that $f_{\lambda}(P)\subseteq{\rm int}(P)$ for all $\lambda\in J$. Without loss of generality we can take $J$ to be a closed interval. For every $\lambda \in J$ we have an attractor
\[
K_\lambda=\bigcap_{n\geq 0}f_\lambda^n(P)
\] with trapping region $P$. We say that the family $\{K_\lambda\}_{\lambda\in J}$ is a \emph{local continuation of attractors} determined by $P$.

In general the topological type of the attractors $K_{\lambda}$ changes with $\lambda$. The main result in this section states that, by contrast, the topological type of the (completions of the) basins of attraction $\hat{\mathcal{A}}(K_{\lambda})$ remains the same. That is, if we fix an arbitrary $\lambda_0 \in J$, we have a family of homeomorphisms \[h_{\lambda} : \hat{\mathcal{A}}(K_{\lambda_0}) \cong \hat{\mathcal{A}}(K_{\lambda}) \quad \text{ for every } \lambda \in J,\] with $h_{\lambda_0} = {\rm id}$. 

We will prove that, in fact, the homeomorphisms $h_{\lambda}$ can be chosen to depend continuously on $\lambda$. To make sense of this last assertion it is convenient to reformulate the setup slightly. Consider the extended phase space $M \times J$ and the map $F(x,\lambda)$ defined there. Notice that this restriction is an open embedding, so we are in the framework of Section \ref{sec:completion}. Clearly $P \times J$ is a trapping region for $F$, and the attractor it determines is just \[K = \bigcap_{n \geq 0} F^n(P \times J) = \bigcup_{\lambda \in J} K_{\lambda} \times \{\lambda\}.\] The basin of attraction of $K$ in $M \times J$ is precisely $\bigcup_{\lambda \in J} \mathcal{A}(K_{\lambda}) \times \{\lambda\}$. Let $\hat{\mathcal{A}}(K)$ be the completion of the basin of attraction of $K$ in $M \times J$. Because $F$ is level-preserving, points in $\hat{\mathcal{A}}(K)$ have a well defined $\lambda$-coordinate, and it should be clear that the $\lambda$-slice of $\hat{\mathcal{A}}(K)$ is the completion of the basin of attraction of $K_{\lambda}$. We will just write $\hat{\mathcal{A}}(K) = \bigcup_{\lambda \in J} \hat{\mathcal{A}}(K_{\lambda}) \times \{\lambda\}$.

With this notation, the main theorem of this section reads as follows:

\begin{theorem}\label{thm:localcont}
      Let $f_\lambda$ be a parametrized family of embeddings of a boundariless manifold $M$. Suppose that $P$ is a trapping region for all the $f_{\lambda}$, $\lambda \in J$. Let $K_{\lambda}$ be the family of attractors determined by $P$ and set $K$ as defined above.  Fix any $\lambda_0 \in J$. There exists a homeomorphism $h:\hat{\mathcal{A}}(K_{\lambda_0}) \times J \longrightarrow \hat{\mathcal{A}}(K)$ such that:
\begin{itemize}
	\item[(i)] for every $\lambda\in J$ 
        \[
        h(\hat{\mathcal{A}}(K_{\lambda_0}) \times\{\lambda\}) =\hat{\mathcal{A}}(K_{\lambda}) \times\{\lambda\},\]
	\item[(ii)] $h$ is the identity on $\hat{\mathcal{A}}(K_{\lambda_0}) \times\{\lambda_0\}$,
	\item[(iii)] $h$ is the identity on $P\times J$. 
\end{itemize}
\end{theorem}

In order to prove Theorem~\ref{thm:localcont} we shall use the following version of a deep result by Edwards and Kirby \cite[Corollary~1.2]{edwardskirby1}.

\begin{theorem}[Edwards-Kirby]\label{thm:edkir} 
Let $g_t: C\longrightarrow W$, $t\in[0,1]$ be an isotopy of a compact subset $C$ of a boundariless manifold $W$ such that $g_0$ is the inclusion of $C$ in $W$. Assume that the isotopy has an extension to a neighbourhood $U$ of $C$. Then there exists an ambient isotopy $\tilde{g}_t:W\longrightarrow W$ such that $\tilde{g}_0={\rm id}$ and $\tilde{g}_t=g_t|_C$. Futhermore, $g_t$ can be chosen to be stationary outside a compact subset of $W$.
\end{theorem}

\begin{proof}[Proof of Theorem~\ref{thm:localcont}] For notational simplicity, and without loss of generality, we take $J = [0,1]$ and $\lambda_0 = 0$. Set $W := {\rm int}_M P$, which is a boundariless manifold. The proof is by induction. We first describe the inductive step in the abstract.

Suppose that $H: P \times [0,1] \longrightarrow P \times [0,1]$ is a level-preserving homeomorphism such that: (i) $h_0 = {\rm Id}_P$, (ii) there is a compact set $B \subseteq W$ such that each $h_{\lambda}$ is the identity outside $B$. As usual, we write $H(x,\lambda) = (h_{\lambda}(x),\lambda)$. Let $C := f_0(P)$ and consider the maps $g_{\lambda} := f_{\lambda} \circ h_{\lambda} \circ f_0^{-1} : C \longrightarrow W$. These give an isotopy of $C$ in $W$ whose initial stage $g_0$ is just the inclusion of $C$ in $W$. Now we check that it has an extension to a neighbourhood $U$ of $C$. Let $V$ be a neighbourhood of $P$ so small that $f_{\lambda}(V) \subseteq W$ for all $\lambda \in [0,1]$. Set $U := f_0(V)$, which is a neighbourhood of $C$ contained $W$. Notice that assumption (ii) on $H$ implies that it can be extended to a level preserving homeomorphism of all of $M \times [0,1]$ by letting it be the identity outside $B \times [0,1]$. Then the same formula that defines $g_{\lambda}$ makes sense on all $U = f_0(V)$ and sends it to $W$. By Theorem \ref{thm:edkir} there exists an isotopy $\tilde{g}_{\lambda}$ of $W$ such that $\tilde{g}_0 = {\rm id}$ and $f_{\lambda} \circ h_{\lambda} \circ f_0^{-1} = g_{\lambda}$ on $C$. Therefore $f \circ h_{\lambda} = \tilde{g}_{\lambda} \circ f_0$ on $P$. If we now set $H' : P \times [0,1] \longrightarrow P \times [0,1]$ to be $H'(x,\lambda) := (\tilde{g}_{\lambda}(x),\lambda)$, we see that $F \circ H = H' \circ (f_0 \times {\rm id})$. Notice also that we can choose the $\tilde{g}_{\lambda}$ to be stationary outside some compact subset of $W$, and so this procedure returns an $H'$ that satisfies the same assumptions listed above for $H$.

Now we can prove the theorem. Let $P_0,P_1,\ldots$ be countably many disjoint copies of $P$. Start with $H^{(0)}$ the identity of $P \times [0,1]$ and apply the above construction inductively to get a family of level-preserving homeomorphisms $H^{(i)}$ of each $P_i \times [0,1]$. These satisfy $F \circ H^{(i)} = H^{(i+1)} \circ (f_0 \times {\rm id})$. Set $Q := \coprod P_i \times [0,1]$. The union of the homeomorphisms $H^{(i)}$ gives a homeomorphism $h : Q \longrightarrow Q$. If one now considers the equivalence relation generated by $(x,\lambda) \in P_i \times [0,1] \sim (f_0(x),\lambda) \in P_{i+1} \times [0,1]$, one sees that $h$ carries it onto the equivalence relation $(x,\lambda) \in P_i \times [0,1] \sim (f_{\lambda}(x),\lambda) \in P_{i+1} \times [0,1]$. Thus $h$ descends to a level-preserving homeomorphism $\bar{h}$ between the quotients of $Q$ by these two relations, which are respectively $\hat{\mathcal{A}}(K_0) \times [0,1]$ and $\hat{\mathcal{A}}(K)$. (The latter identification follows from Remark \ref{rem:completionP}). The condition that all the $H^{(i)}$ be the identity at level zero ensures that $\bar{h}$ is the identity on $\hat{\mathcal{A}}(K_0)$. Finally, recall that there is a homeomorphic copy of $P$ in the completion of the basins of attraction represented by points of the form $[x \in P_0]$. Since $h$ is $H^{(0)} = {\rm Id}$ on $P_0 \times [0,1]$, it follows that $\bar{h}$ is the identity on $P \times [0,1]$. 
\end{proof}

The previous local result can be extended to any continuous family of attractors of embeddings. To make this precise, let $\{K_\lambda\}_{\lambda\in[0,1]}$ be a family of subcompacta of $M$ such that for each $\lambda$, $K_\lambda$ is an attractor of $f_\lambda$. We say that the family $\{K_\lambda\}_{\lambda\in[0,1]}$ is a \emph{continuation of attractors} if for every $\lambda_0\in[0,1]$ there exist a trapping region $P$ of $K_{\lambda_0}$ and a non-degenerate subinterval $J$ with $\lambda_0\in{\rm int}_{[0,1]} J$ such that the subfamily $\{K_\lambda\}_{\lambda\in J}$ is the local continuation of $K_{\lambda_0}$ determined by $P$. Mimicking the construction described above before Theorem \ref{thm:localcont}, but this time performed locally on each $J$, one sees that $K = \bigcup_{\lambda \in [0,1]} K_{\lambda}$ is an attractor for the map $F$ on $M \times [0,1]$ (see the proof of the next corollary for more details).

\begin{corollary}\label{coro:globalcont}
      In the setting just described, there exists a homeomorphism $h:\hat{\mathcal{A}}(K_0)\times[0,1]\longrightarrow \hat{\mathcal{A}}(K)$ such that for every $\lambda\in[0,1]$,
        \[
        h(\hat{\mathcal{A}}(K_0)\times\{\lambda\})=\hat{\mathcal{A}}(K_\lambda) \times\{\lambda\}.
        \]
\end{corollary}

\begin{proof} 
For the sake of brevity we denote by $\hat{\mathcal{A}}_{\lambda}$ the completion of the basin of attraction of each $K_{\lambda}$, and similarly set $\hat{\mathcal{A}}$ the completion of the basin of attraction of $K$. Since $[0,1]$ is compact and connected there exists a finite subdivision $0=\lambda_0<\lambda_1<\lambda_2<\ldots<\lambda_k=1$ of $[0,1]$ and compact sets $P_0,P_1,\ldots,P_{k-1}$ such that $P_i$ is a trapping region that determines the local continuation $\{K_\lambda\}_{\lambda\in[\lambda_{i},\lambda_{i+1}]}$. It is not difficult to see that compact set
\[
P=\bigcup_{i=0}^{k-1} P_i\times[\lambda_i,\lambda_{i+1}]
\]
is a trapping region for $K$ and, hence, $K$ is an attractor for $F$. We can then apply Theorem~\ref{thm:localcont} to each piece $J = [\lambda_i,\lambda_{i+1}]$ to obtain homeomorphisms $h_i:\hat{\mathcal{A}}_{\lambda_i}\times[\lambda_i,\lambda_{i+1}]\longrightarrow\hat{\mathcal{A}}^i$, where $\hat{\mathcal{A}^i}$ is the completion of the basin of attraction of $K \cap M \times [\lambda_i,\lambda_{i+1}]$. Consider the map $\widetilde{h}_i(x)=p_1\circ h_i|_{\hat{\mathcal{A}}_{\lambda_i}\times\{\lambda_{i+1}\}}(x,\lambda_{i+1}) :\hat{\mathcal{A}}_{\lambda_i}\longrightarrow \hat{\mathcal{A}}_{\lambda_{i+1}}$ where $p_1:M\times[0,1]\longrightarrow M$ denotes the projection onto the first factor. By construction $\widetilde{h}_i$ is a homeomorphism and
\begin{equation}
h_i(x,\lambda_{i+1})=(\widetilde{h}_i(x),\lambda_{i+1})
\label{eq:cont}
\end{equation}

We define the homeomorphism $h:\hat{\mathcal{A}}_0\times[0,1]\longrightarrow \hat{\mathcal{A}}$ in each piece $\hat{\mathcal{A}}_0\times[\lambda_i,\lambda_{i+1}]$ as
\[
h(x,\lambda)=h_i(\widetilde{h}_{i-1}\circ\ldots\circ\widetilde{h}_0(x),\lambda).
\]
This map is a well defined homeomorphism by \eqref{eq:cont} together with the fact that $h_i$ is the identity in $\hat{\mathcal{A}}_{\lambda_i}\times\{\lambda_i\}$.

\end{proof}

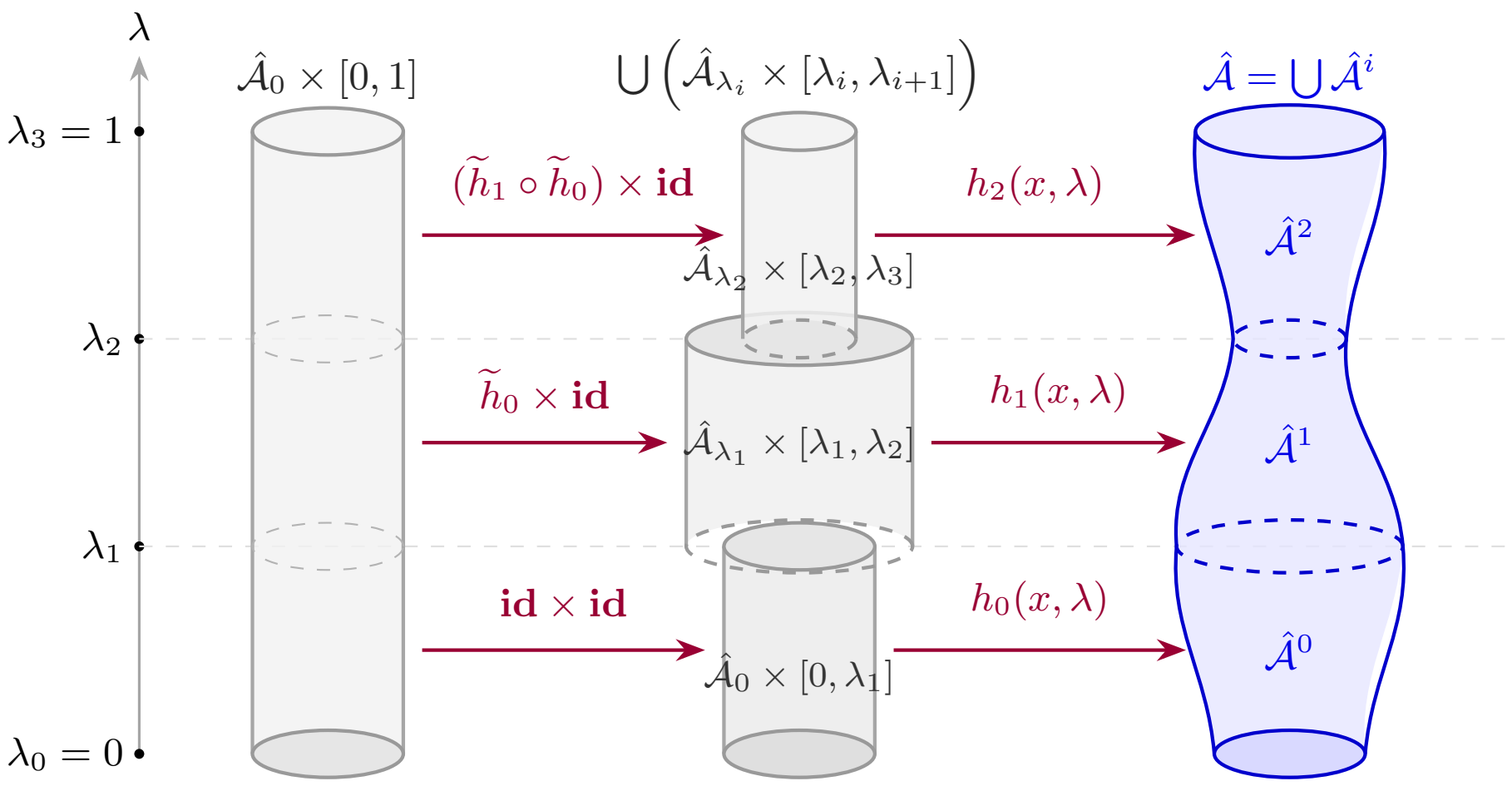
\begin{figure}[htbp]
\centering
\begin{tikzpicture}[scale=0.75, >=Stealth]

  \draw[->, gray!70, semithick] (-1.8, 0) -- (-1.8, 7.4) node[above, black, font=\small] {$\lambda$};

  \fill[black] (-1.8, 0.0) circle (1.5pt) node[left, font=\small] {$\lambda_0=0$};
  \fill[black] (-1.8, 2.2) circle (1.5pt) node[left, font=\small] {$\lambda_1$};
  \fill[black] (-1.8, 4.4) circle (1.5pt) node[left, font=\small] {$\lambda_2$};
  \fill[black] (-1.8, 6.6) circle (1.5pt) node[left, font=\small] {$\lambda_3=1$};

  \draw[gray!25, dashed] (-1.8, 2.2) -- (12.8, 2.2);
  \draw[gray!25, dashed] (-1.8, 4.4) -- (12.8, 4.4);

  \begin{scope}[xshift=0.2cm]

    \fill[gray!10, opacity=0.8] 
      (-0.8, 0) -- (-0.8, 6.6) arc[start angle=180, end angle=0, x radius=0.8, y radius=0.25] -- (0.8, 0) arc[start angle=0, end angle=180, x radius=0.8, y radius=0.25] -- cycle;

    \draw[thick, gray!70] (-0.8, 0) -- (-0.8, 6.6);
    \draw[thick, gray!70] (0.8, 0) -- (0.8, 6.6);

    \draw[thick, gray!80, fill=gray!20] (0, 0) ellipse (0.8cm and 0.25cm);
    \draw[gray!60, dashed] (0, 2.2) ellipse (0.8cm and 0.25cm);
    \draw[gray!60, dashed] (0, 4.4) ellipse (0.8cm and 0.25cm);
    \draw[thick, gray!80] (0, 6.6) ellipse (0.8cm and 0.25cm);

    \node[font=\small\bfseries, black!85] at (0, 7.2) {$\hat{\mathcal{A}}_0 \times [0,1]$};

  \end{scope}

  \draw[->, thick, purple!80!black] (1.2, 1.1) -- (4.2, 1.1) 
    node[midway, above=3pt, font=\small\bfseries] {$\text{id} \times \text{id}$};

  \draw[->, thick, purple!80!black] (1.2, 3.3) -- (3.8, 3.3) 
    node[midway, above=3pt, font=\small\bfseries] {$\widetilde{h}_0 \times \text{id}$};

  \draw[->, thick, purple!80!black] (1.2, 5.5) -- (4.4, 5.5) 
    node[midway, above=3pt, font=\small\bfseries] {$(\widetilde{h}_1 \circ \widetilde{h}_0) \times \text{id}$};

  \begin{scope}[xshift=5.2cm]

    \fill[gray!15, opacity=0.85] 
      (-0.8, 0) -- (-0.8, 2.2) arc[start angle=180, end angle=0, x radius=0.8, y radius=0.25] -- (0.8, 0) arc[start angle=0, end angle=180, x radius=0.8, y radius=0.25] -- cycle;
    \draw[thick, gray!70] (-0.8, 0) -- (-0.8, 2.2);
    \draw[thick, gray!70] (0.8, 0) -- (0.8, 2.2);
    \draw[thick, gray!80, fill=gray!25] (0, 0) ellipse (0.8cm and 0.25cm);
    \draw[thick, gray!80, fill=gray!20] (0, 2.2) ellipse (0.8cm and 0.25cm);
    \node[font=\footnotesize\bfseries, black!80] at (0, 0.85) {$\hat{\mathcal{A}}_0 \times [0,\lambda_1]$};

    \fill[gray!12, opacity=0.85] 
      (-1.2, 2.2) -- (-1.2, 4.4) arc[start angle=180, end angle=0, x radius=1.2, y radius=0.28] -- (1.2, 2.2) arc[start angle=0, end angle=180, x radius=1.2, y radius=0.28] -- cycle;
    \draw[thick, gray!70] (-1.2, 2.2) -- (-1.2, 4.4);
    \draw[thick, gray!70] (1.2, 2.2) -- (1.2, 4.4);
    \draw[thick, gray!80, dashed] (0, 2.2) ellipse (1.2cm and 0.28cm);
    \draw[thick, gray!80, fill=gray!20] (0, 4.4) ellipse (1.2cm and 0.28cm);
    \node[font=\footnotesize\bfseries, black!80] at (0, 3.3) {$\hat{\mathcal{A}}_{\lambda_1} \times [\lambda_1,\lambda_2]$};

    \fill[gray!10, opacity=0.85] 
      (-0.6, 4.4) -- (-0.6, 6.6) arc[start angle=180, end angle=0, x radius=0.6, y radius=0.20] -- (0.6, 4.4) arc[start angle=0, end angle=180, x radius=0.6, y radius=0.20] -- cycle;
    \draw[thick, gray!70] (-0.6, 4.4) -- (-0.6, 6.6);
    \draw[thick, gray!70] (0.6, 4.4) -- (0.6, 6.6);
    \draw[thick, gray!80, dashed] (0, 4.4) ellipse (0.6cm and 0.20cm);
    \draw[thick, gray!80] (0, 6.6) ellipse (0.6cm and 0.20cm);
    \node[font=\footnotesize\bfseries, black!80] at (0, 5.15) {$\hat{\mathcal{A}}_{\lambda_2} \times [\lambda_2,\lambda_3]$};

    \node[font=\small\bfseries, black!85] at (0, 7.2) {$\bigcup \left(\hat{\mathcal{A}}_{\lambda_i} \times [\lambda_i, \lambda_{i+1}]\right)$};

  \end{scope}

  \draw[->, thick, purple!80!black] (6.2, 1.1) -- (9.3, 1.1) 
    node[midway, above=3pt, font=\small\bfseries] {$h_0(x, \lambda)$};

  \draw[->, thick, purple!80!black] (6.6, 3.3) -- (9.3, 3.3) 
    node[midway, above=3pt, font=\small\bfseries] {$h_1(x, \lambda)$};

  \draw[->, thick, purple!80!black] (6.0, 5.5) -- (9.4, 5.5) 
    node[midway, above=3pt, font=\small\bfseries] {$h_2(x, \lambda)$};

  \begin{scope}[xshift=10.4cm]

    \fill[blue!10, opacity=0.8] 
      (-0.8, 0) to[out=95, in=265] (-1.2, 2.2) to[out=95, in=265] (-0.6, 4.4) to[out=95, in=265] (-1.0, 6.6)
      arc[start angle=180, end angle=0, x radius=1.0, y radius=0.28]
      to[out=265, in=95] (0.6, 4.4) to[out=265, in=95] (1.2, 2.2) to[out=265, in=95] (0.8, 0)
      arc[start angle=0, end angle=180, x radius=0.8, y radius=0.25] -- cycle;

    \draw[thick, blue!80!black] 
      (-0.8, 0) to[out=95, in=265] (-1.2, 2.2) to[out=95, in=265] (-0.6, 4.4) to[out=95, in=265] (-1.0, 6.6);
    \draw[thick, blue!80!black] 
      (0.8, 0) to[out=85, in=275] (1.2, 2.2) to[out=95, in=265] (0.6, 4.4) to[out=85, in=275] (1.0, 6.6);

    \draw[thick, blue!80!black, fill=blue!15] (0, 0) ellipse (0.8cm and 0.25cm);
    \draw[thick, blue!80!black, dashed] (0, 2.2) ellipse (1.2cm and 0.28cm);
    \draw[thick, blue!80!black, dashed] (0, 4.4) ellipse (0.6cm and 0.20cm);
    \draw[thick, blue!80!black] (0, 6.6) ellipse (1.0cm and 0.28cm);

    \node[font=\small\bfseries, blue!90!black] at (0, 7.2) {$\hat{\mathcal{A}} = \bigcup \hat{\mathcal{A}}^i$};
    \node[blue!90!black, font=\small\bfseries] at (0, 1.1) {$\hat{\mathcal{A}}^0$};
    \node[blue!90!black, font=\small\bfseries] at (0, 3.3) {$\hat{\mathcal{A}}^1$};
    \node[blue!90!black, font=\small\bfseries] at (0, 5.5) {$\hat{\mathcal{A}}^2$};

  \end{scope}

\end{tikzpicture}
\caption{Construction of the basin homeomorphism $h: \hat{\mathcal{A}}_0 \times [0,1] \longrightarrow \hat{\mathcal{A}}$ from Corollary~\ref{coro:globalcont}}
\label{fig:continuation_two_steps_id}
\end{figure}

We conclude this section with a comment. An attractor is an isolated set in the sense of Conley, and its homological Conley index actually coincides with the homology of the completed basin of attraction. Since the Conley index is invariant under continuations, it follows that $H_*(\hat{\mathcal{A}}_{\lambda})$ is also invariant under continuations. The theorem above is stronger because it states that the topological type of the $\hat{\mathcal{A}}_{\lambda}$ is invariant.

\section{The basics of generalized Hopf bifurcations} \label{sec:hopf0}

The Poincaré-Andronov-Hopf and the Neimark-Sacker bifurcations are the continuous and discrete-time versions of a dynamical phenomenon, which, in its most basic form, goes as follows. One has a parameterized family of dynamical systems (ODEs or diffeomorphisms) in $\mathbb{R}^2$ which all have a point $p$ as an equilibrium. For a certain value of the parameter $\lambda = \lambda_0$ there is a change in the stability of $p$: for $\lambda < \lambda_0$ the point is a hyperbolic attractor and for $\lambda > \lambda_0$ it becomes a hyperbolic repeller. Under suitable assumptions on how this change in stability takes place, the classical theory predicts that a family of attracting, closed invariant curves will grow out (``be expelled by'') the fixed point $p$ as $\lambda$ grows past $\lambda_0$. One of those assumptions concerns the normal form of the dynamics at $\lambda = \lambda_0$ and implies that $p$ is also an attractor for that value of the parameter (although not a hyperbolic one). Thus, if we focus only on the change in stability of $p$ and disregard how this change happens, the situation may be summarized by saying that the point $p$ is an attractor for $\lambda \leq \lambda_0$ and becomes a repeller for $\lambda > \lambda_0$. From this perspective it makes sense to replace $p$ with an arbitrary compact set $K$, therefore allowing completely general objects to undergo this kind of bifurcation. This was called by Seiberg and Florio \cite{MR1277591} a generalized Poincaré-Andronov-Hopf bifurcation. We will just speak of generalized Hopf bifurcations for short.

The definition proposed by Seibert and Florio makes sense also for semidynamical systems, and it is worth studying it in this generality. The reason is that then one can include in the analysis situations where the change in stability occurs only inside some invariant subspace of ambient space, or more generally inside some parameter-dependent locally invariant subspace of ambient space. To explain this, consider first the following very simple example:

\begin{example} \label{ex:motivation1} Suppose $f_{\lambda}$ is a continuous family of $\mathcal{C}^1$ diffeomorphisms of $\mathbb{R}^3$, all having a point $p$ as a fixed point. Assume that for $\lambda < \lambda_0$ the derivative $d_p f_{\lambda}$ has one real and two complex conjugate eigenvalues, all with a modulus less than $1$. Suppose also that as $\lambda$ increases the real eigenvalue remains always less than $1$ in modulus, but the complex conjugate pair increases its modulus, reaching $1$ at $\lambda = \lambda_0$ and growing larger than $1$ afterwards. In this situation, for $\lambda \neq \lambda_0$ the point $p$ is hyperbolic: it is an attractor for $\lambda < \lambda_0$, whereas for $\lambda > \lambda_0$ it still attracts along a one-dimensional (local) invariant manifold $W^s_{\lambda}$ but now repels along a two-dimensional (local) invariant manifold $W^u_{\lambda}$.
\end{example}

This example does not fit the above definition of a generalized Hopf bifurcation because $p$ is not a repeller in $\mathbb{R}^3$ for $\lambda > \lambda_0$. However, we can accomodate it in that
framework provided we do not insist on working with flows or homeomorphisms but expand our scope slightly to injective semidynamical systems. Let us illustrate how with the same example:

\begin{example}[Example \ref{ex:motivation1} continued] By applying the center manifold theorem in the extended phase space (i.e. $\mathbb{R}^3$ enlarged with a $\lambda$-axis) one obtains a continuous family of two-dimensional local invariant manifolds $\lambda \longmapsto W_{\lambda}$ all of which contain $p$ and such that $p$ is an attractor inside $W_{\lambda}$ for $\lambda < \lambda_0$ and a repeller in $W_{\lambda}$ for $\lambda > \lambda_0$. As in the classical case mentioned earlier, appropriate higher order assumptions will ensure that $p$ is an attractor inside $W_{\lambda}$ for $\lambda = \lambda_0$ as well. Now, the assertion that $\lambda \longmapsto W_{\lambda}$ is a continuous family means that there exist an open $2$-dimensional manifold $W$ and a continuous mapping $w : W \times (\lambda_0 - \epsilon,\lambda_0 + \epsilon) \longrightarrow \mathbb{R}^3$ which, for every $\lambda \in (\lambda_0-\epsilon,\lambda_0+\epsilon)$, is a homeomorphism onto $W_{\lambda}$. Let us take $W$ to be the unit disk in the plane for definiteness, and choose $w(\lambda,0) = p$. Also, let us denote by $w_{\lambda} : W \longrightarrow W_{\lambda}$ the map $x \longmapsto w(\lambda,x)$. For $\lambda = \lambda_0$ the point $p$ is an attractor in $W_{\lambda_0}$ and so we can find a small closed disk $D \subseteq W$ centered at $0$ and such that $f_{\lambda_0}(w_{\lambda_0}(D)) \subseteq w_{\lambda_0}({\rm int}\ D)$. This condition will still hold true for $\lambda$ sufficiently close to $\lambda_0$, and so for such values of the parameter we can consider the maps $g_{\lambda} : D \longrightarrow D$ given by $g_{\lambda}:= w_{\lambda}^{-1} \circ f_{\lambda} \circ w_{\lambda}$. These are all embeddings, but obviously not onto. The point $0$ is an attractor for $g_{\lambda}$ with $\lambda \leq \lambda_0$ and a repeller for $\lambda > \lambda_0$. Any invariant sets that bifurcate from $p$ under $f_{\lambda}$ will lie in $w_{\lambda_0}(D)$ for small enough $\lambda$ and therefore manifest themselves as invariant sets that bifurcate from $0$ in $D$ under $g_{\lambda}$ (and conversely). In other words, the study of the bifurcation properties of $p$ for the diffeomorphisms $f_{\lambda}$ can be reduced to the study of a bifurcation in the boundariless manifold $D$ for the embeddings (not diffeomorphisms any longer) $g_{\lambda}$.
\end{example}

Summing up, and with a trivial rescaling of the parameter $\lambda$, our definitive setup is as follows. The ambient space is a boundariless manifold $M$ and $f_{\lambda} : M \longrightarrow M$ is a continuous family of embeddings, with $\lambda \in [0,1]$. There is a compact set $K \subseteq M$ that is invariant for all the $f_{\lambda}$, is an attractor for $\lambda = 0$ and becomes a repeller for $\lambda > 0$. We will then say that $K$ undergoes a \emph{generalized Hopf bifurcation}.

We emphasize that we take this definition to sit in the realm of Topological Dynamics, so that ``attractor'' is meant in the sense of Section \ref{sec:completion}. Similarly for ``repeller'', although since the dynamics may not be invertible a few comments are in order. A repeller is the analogue of an attractor but in reverse time.  Suppose $f$ is an embedding and $N$ is a compact set such that $N \subseteq f({\rm int}\ N)$. Then every point in $N$ has a preimage in ${\rm int}\ N$; it must be unique because $f$ is an embedding. Hence the map $(f|_N)^{-1}$ is well defined and provides a homeomorphic embedding of $N$ into its interior. There is an attractor $K$ in $N$ for the time-reversed dynamics $(f|_N)^{-1}$, namely $K:= \bigcap_{n \geq 0} (f|_N)^{-n}(N)$. This is what we call a repeller of the forward-time dynamics $f$. Its basin of repulsion $\mathcal{R}(K)$ is the set of points which have a full backward orbit that converges to $K$ asymptotically. Bearing in mind that $N \subseteq f(N)$ implies that points in $N$ do have a full backwards orbit, it is straightforward to check that the basin of repulsion is $\mathcal{R}(K) = \bigcup_{n \geq 0} f^n(N)$.

Now suppose again that $f_{\lambda}$ is a continuous family of embeddings and $K$ is a compact set that undergoes a generalized Hopf bifurcation at $\lambda = 0$, as defined above. Let $P$ be a trapping region for $f_0$ that determines the attractor $K$; i.e. $f_0(P) \subseteq {\rm int}(P)$. We have that $P$ is still a trapping region for all $\lambda \approx 0$, and so each $f_{\lambda}$ determines an attractor $K_{\lambda}$ inside $P$. We claim that the basin of repulsion $\mathcal{R}_{\lambda}(K)$ is contained in $K_{\lambda}$. Recall that any point $x \in \mathcal{R}_{\lambda}(K)$ has a full negative semiorbit that converges to $K$, so in particular it eventually enters $P$. Thus $x \in \bigcap f_{\lambda}^n(P) = K_{\lambda}$. Set $A_{\lambda} := K_{\lambda} \setminus \mathcal{R}_{\lambda}(K)$, which is an attractor in $K_{\lambda}$. Since $K_{\lambda}$ is an attractor in ambient space, $A_{\lambda}$ is also an attractor in ambient space \cite[Proposition~5]{MR4421505}. Obviously the $A_{\lambda}$ are all disjoint from $K$. Also, an easy argument shows that they converge to $K$ as $\lambda \rightarrow 0$ in the sense that for any neighbourhood $U$ of $K$ there is some $\lambda_U$ such that $A_{\lambda} \subseteq U$ for every $0 < \lambda < \lambda_U$. This is appropriately summarized by saying that the $A_{\lambda}$ are the attractors ``expelled'' by $K$ during this bifurcation. At least two natural questions arise: what is the nature of the expelled $A_{\lambda}$?, and can any $K$ undergo such a bifurcation?

Before moving ahead we show that generalized Hopf bifurcations can indeed happen beyond the classical cases. In fact, the next proposition shows that any attractor for a flow can undergo a generalized Hopf bifurcation. The proof of the proposition helps gain some intuition on how these can look like. A very interesting fact is that the analogue for discrete dynamics is not true: there are attractors for homeomorphisms that cannot undergo Hopf bifurcations. We will discuss this in detail in Section \ref{sec:compacta_that_cannot}.

\begin{proposition} \label{prop:hopfad} Suppose $K$ is an attractor for an injective semiflow $\varphi$ on $M$. Then there exists a family of injective semiflows $\varphi_{\lambda}$, $\lambda \in [0,1]$, such that $K$ undergoes a generalized Hopf bifurcation for $\varphi_{\lambda}$.
\end{proposition}
\begin{proof}
We shall prove the theorem assuming that $\varphi$ is a flow. The case when $\varphi$ is a semiflow can be reduced to this by using the completion of the basin of attraction described in Section~\ref{sec:semi}.

Since $K$ is a stable attractor for $\varphi$, it has a proper Lyapunov function $L : \mathcal{A}(K) \longrightarrow [0,+\infty]$. Recall that $L$ has this properties: (i) it is continuous, (ii) $L$ is strictly decreasing along the trajectories in $\mathcal{A}(K) \setminus K$, (iii) $L$ takes the value $0$ precisely on $K$, and (iv) $L^{-1}([0,c])$ is compact for every $c$. Notice that property (ii) implies that any level set $L^{-1}(c)$, $c > 0$, is a section of the flow in $\mathcal{A}(K) \setminus K$: every trajectory there intersects exactly once that set.

Extend $\varphi$ to a level-preserving flow $\Phi$ on $M \times [0,1]$ in the obvious way $\Phi((x,u),\lambda) := (\varphi(x,t), \lambda)$. This has $K \times [0,1]$ as an attractor with basin of attraction $\mathcal{A}(K) \times [0,1]$. Define the set \[S := \{(x,\lambda) \in \mathcal{A}(K) \times [0,1] : L(x) = \lambda\}.\] This is closed in $M \times [0,1]$ and, in fact, compact by property (iv). Now we modify the flow $\Phi$ by slowing it down in such a manner that it becomes stationary on $S$. This is very easy to do for a flow generated by an ODE $\dot{x} = f(x)$; it suffices to multiply $f(x)$ by a function $\theta(x)$ that is nonnegative and zero exactly on $S$. In the purely topological context it can also be done but the argument is more complicated and can be found in \cite[Chapters 4 and 5]{beck2}. The new flow $\Psi$ on $M \times [0,1]$ is still level-preserving, and we denote by $\psi_{\lambda}$ the flow on $M$ given by the $\lambda$-slice of $\Psi$. By construction, every point in the set $L^{-1}(\lambda)$ is stationary for $\psi_{\lambda}$. The following assertions should be more or less obvious:

(1) For $\psi_0$ every point in $K$ is stationary and the trajectory of every other point is just an orientation-preserving reparameterization of its original trajectory under $\varphi$. In particular $L$ is still a Lyapunov function, and so $K$ is still an attractor for $\psi_0$.

(2) For $0 < \lambda \leq 1$ the trajectories of $\psi_{\lambda}$ are (orientation preserving reparameterizations of) the trajectories of $\varphi$ with the exception of the trajectories in $\mathcal{A}(K) \setminus K$, which get split into two at their unique point of intersection with $S$. In particular $L$ is still strictly decreasing along the trajectories of $\psi_{\lambda}$ in $\mathcal{A}(K) \setminus K$ other than the stationary points in $S$.

(3) For $0 < \lambda \leq 1$ the set $K$ is an attractor for $\psi_{\lambda}$ as shown by the Lyapunov function $L_{\lambda} : L^{-1}([0,\lambda)) \longrightarrow [0,+\infty)$ defined by \[L_{\lambda}(x) := \left\{ \begin{array}{ll} \frac{1}{\lambda - L(x)} - \frac{1}{\lambda} & \text{ if } x \not\in K \\ 0 & \text{ if } x \in K \end{array} \right.\] The basin of attraction of $K$ is precisely $L^{-1}([0,\lambda))$. 

Finally, we define the flows $\varphi_{\lambda}$ on $M$ by \[\varphi_{\lambda}(x,t) := \left\{ \begin{array}{ll} \psi_{\lambda}(x,-t) & \text{ if } x \in L^{-1}([0,\lambda]) \\ \psi_{\lambda}(x,t) & \text{ otherwise} \end{array} \right.\] that is, we just reverse the dynamics in the basin of attraction of $K$. The two definitions agree on the intersection of their domains $L^{-1}(\lambda)$ where both are stationary, so the $\varphi_{\lambda}$ are continuous. The set $K$ is an attractor for $\varphi_0$ by (1) and becomes a repeller for $\lambda > 0$ by (3). Finally, the family $\varphi_{\lambda}$ depends continuously on $\lambda$ because they are the slices of the flow on $M \times [0,1]$ given by \[((x,\lambda),t) \longmapsto \left\{ \begin{array}{ll} \Psi((x,\lambda),-t) & \text{ if } x \in L^{-1}([0,\lambda]) \\ \Psi((x,\lambda),t) & \text{ otherwise} \end{array} \right.\]
\end{proof}

Notice that the proposition claims that $K$ admits a Hopf bifurcation for some family of flows $\varphi_{\lambda}$ but does not claim that $\varphi_0$ coincides with the original flow $\varphi$. Indeed, the proof of the proposition requires that $\varphi$ be slowed down to a halt on $K$ (at $\lambda = 0$) before $K$ becomes a repeller for $\lambda > 0$. The flow $\varphi_0$ does coincide with $\varphi$, at least qualitatively, outside of $K$. It is very easy to modify the argument of the proof so that $\varphi_0 = \varphi$, in the interval $\lambda \in [0,1/2]$ the flow $\varphi$ gets gradually slowed down to a halt in $K$, and then the generalized Hopf bifurcation takes place at $\lambda = 1/2$, with $K$ becoming a repeller for $\lambda \in (1/2,1]$.

\section{Homological analysis of generalized Hopf bifurcations} \label{sec:homological}

Our goal in this section is to give a homological description of the attractors $A_{\lambda}$ expelled by a compactum $K$ that undergoes a Hopf bifurcation. The spirit of this section is that the cohomology of $K$ is known and we want to find the cohomology of the $A_{\lambda}$ just from that information. We point out that the results of this section are fully general, removing some additional conditions about concentricity and topological spines that are used in \cite{MR4421505} for technical reasons. These conditions turn out to be redundant in manifolds of dimension at most three but seem too strong in higher dimensions. We shall see the redundancy of the topological spine condition in the next section.

\begin{remark} Throught this section we assume that the ambient manifold $M$ is smooth or triangulable, so that we can freely use Theorem \ref{thm:pacojaime}. Of course, no such assumption is placed on the dynamics itself.
\end{remark}

We begin with an auxiliary result that is, nevertheless, interesting in its own right. We will look at it from a different perspective in the following section.

\begin{theorem}\label{thm:finitehopf}
 If $K$ admits a Hopf bifurcation its integral \v Cech cohomology groups are finitely generated and only a finite number of them are nonzero.
\end{theorem}
\begin{proof} Assume that $K$ admits a Hopf bifurcation and that $f_\lambda$ is a family that realizes it. Choose a trapping region $P$ for $f_0$. Now, fix $\lambda>0$ and $n \geq 1$ sufficiently large so that $f_0^n(P)\subseteq\mathcal{R}_\lambda(K)$. For any $k \geq 0$ define $Q_k := f_{\lambda}^{-n-k} \circ f_0^n(P)$. These form a basis of compact neighbourhoods of $K$.

We are going to show that the inclusion induced map $\check{H}^*(P) \longrightarrow \check{H}^*(K)$ is surjective. Let $u \in \check{H}^*(K)$. By the continuity property of \v{C}ech cohomology $u$ admits an extension to a neighbourhood of $K$, which in turn contains some $Q_k$ for large $k$. That is, there exist $k$ and $\hat{u} \in \check{H}^*(Q_k)$ such that $\hat{u}|_K = u$. Consider the commutative diagram \[\xymatrix{P \ar[r]^-{f_0^{n+k}} & f_0^{n+k}(P) \subseteq f_0^n(P) \ar[r]^-{f_{\lambda}^{-n-k}} & Q_k \\ K \ar[u] \ar[rr]^{f_{\lambda}^{-n-k} |_K \circ f_0^{n+k}|_K} & & K \ar[u]}\] The vertical arrows are inclusions. The lower horizontal arrow is homotopic to the identity because the continuation restricts to an isotopy $f_0|_K \simeq f_{\lambda}|_K$. Denote by $g : P \longrightarrow Q_k$ the composition of the maps in the upper row. Then, passing to cohomology the commutativity of the diagram implies that $g^*(z)|_K = z|_K$ for every $z \in \check{H}^*(Q_k)$. In particular, choosing $z = \hat{u}$ we find an element $g^*(\hat{u}) \in H^*(P)$ whose restriction to $K$ is precisely $\hat{u}|_K = u$. Thus the inclusion induced map $\check{H}^*(P) \longrightarrow \check{H}^*(K)$ is surjective.

Since we have assumed $M$ to be a smooth or triangulable manifold, any compact set has a neighbourhood basis of compact manifolds. In particular, the attractor $K$ has arbitrarily small compact neighbourhoods with a finitely generated cohomology. Let $N \subseteq P$ be such a neighbourhood of $K$. Then the inclusion induced map $i^* : \check{H}^*(P) \longrightarrow H^*(K)$ factors through $\check{H}^*(N)$, which coincides with $H^*(N)$. The latter is finitely generated and nonzero only in finitely many degrees, and since $i^*$ is surjective, the same holds for $\check{H}^*(K)$.
  \end{proof}

\begin{theorem}\label{thm:hopfbif}
    Suppose that the attractor $K$ undergoes a Hopf bifurcation. Then, there exists $\lambda_0>0$ such that for $0<\lambda<\lambda_0$ the cohomology of the expelled attractors $A_{\lambda}$ is given by
    \begin{equation} 
    \check{H}^*(A_\lambda)\cong H^*(\hat{\mathcal{A}}\setminus K),
    \label{thm:cohattexp}
    \end{equation}
    where $\hat{\mathcal{A}}$ is the completion of the basin of attraction of $K$. Moreover, these cohomology groups are all finitely generated.
\end{theorem}

We will need the following topological lemma, which we prove after the theorem:

\begin{lemma} \label{lem:recubridor} Let $U$ be a manifold without boundary. Assume $U$ is nonorientable and let $\tilde{U}$ be its orientation cover. If the homology $H_*(\tilde{U})$ is finitely generated, then so is the homology $H_*(U)$.
\end{lemma}

\begin{proof}[Proof of Theorem \ref{thm:hopfbif}] We depart momentarily from the notation of the statement of the theorem and denote as usual $\hat{\mathcal{A}}(K_{\lambda})$ the completion of the basin of attraction of the $K_{\lambda}$, endowed with the homeomorphism $\hat{f}_{\lambda}$. By Theorem \ref{thm:localcont} we have a continuous family of homeomorphisms $h_{\lambda} : \hat{\mathcal{A}}(K_0) \longrightarrow \hat{\mathcal{A}}(K_{\lambda})$. We may then define the continuous family of homeomorphisms $\hat{g}_{\lambda} := h_{\lambda}^{-1} \circ \hat{f}_{\lambda} \circ h_{\lambda}$ on $\hat{\mathcal{A}}(K_0)$. Recall that $h_{\lambda}$ is the identity on $P$, which is a neighbourhood of $K$ that also contains the expelled attractors $A_{\lambda}$ for $\lambda > \lambda_0$ sufficiently close to $\lambda_0$. It follows that $K$, regarded now as a subset of $\hat{\mathcal{A}}(K_0)$, also undergoes a Hopf bifurcation for the family $\hat{g}_{\lambda}$ and the attractors it expels are exactly the same $A_{\lambda}$. In other words, to prove the theorem we can (and will) assume that the Hopf bifurcation takes place in $\hat{\mathcal{A}}(K_0)$ and the dynamics is given by a family of homeomorphisms (the $\hat{g}_{\lambda}$). 

We now revert to the simpler notation of the statement of the theorem: $K_0 = K$ and $\hat{\mathcal{A}}(K_0) = \hat{\mathcal{A}}$. Notice that $A_{\lambda}$ is an attractor for $\hat{g}_{\lambda}$ whose basin of attraction is $\hat{\mathcal{A}} \setminus K$. Thus by Theorem \ref{thm:pacojaime} to prove the theorem it suffices to show that $\hat{\mathcal{A}} \setminus K$ has a finitely generated cohomology. Since $K$ has a finitely generated cohomology by Theorem \ref{thm:finitehopf}, again Theorem \ref{thm:pacojaime} implies that $\hat{\mathcal{A}}$ has a finitely generated cohomology. From the long exact sequence for the pair $(\hat{\mathcal{A}},\hat{\mathcal{A}} \setminus K)$ it follows that it is enough to show that $H^*(\hat{\mathcal{A}},\hat{\mathcal{A}} \setminus K)$ is finitely generated.

Recall that $\hat{\mathcal{A}}$ is an $n$-manifold, and suppose first that it is orientable. Then by Alexander duality $H_q(\hat{\mathcal{A}},\hat{\mathcal{A}} \setminus K) \cong \check{H}^{n-q}(K)$ which by Theorem \ref{thm:finitehopf} is finitely generated. Thus in this case we are finished. For the nonorientable case we reason as follows. Let $p : \tilde{\hat{\mathcal{A}}} \longrightarrow \hat{\mathcal{A}}$ be the orientation cover of $\hat{\mathcal{A}}$ and set $\tilde{K} := p^{-1}(K)$. The family of homeomorphisms $\hat{g}_{\lambda}$ lifts to a family of homeomorphisms $\tilde{g}_{\lambda}$ which give a Hopf bifurcation of $\tilde{K}$ in $\tilde{\hat{\mathcal{A}}}$. Reasoning as in the previous paragraph we deduce that $\tilde{\hat{\mathcal{A}}} \setminus \tilde{K}$ has a finitely generated cohomology. This is a double cover of $\hat{\mathcal{A}} \setminus K$, so the latter also has a finitely generated cohomology by Lemma \ref{lem:recubridor}. This finishes the proof.
\end{proof}

\begin{proof}[Proof of Lemma \ref{lem:recubridor}] Let $p : \tilde{U} \longrightarrow U$ be the covering projection. It is a 2-to-1 map which induces $p_* : H_*(\tilde{U}) \longrightarrow H_*(U)$. In the opposite direction we have the so-called transfer map (see \cite[\S 3.G]{hatcher1}). Roughly, for each singular simplex in $U$ one considers the sum of its two lifts in $\tilde{U}$. This induces a chain map and eventually a map in homology $\pi_* : H_*(U) \longrightarrow H_*(\tilde{U})$. The composition $p_* \circ \pi_*$ is obviously multiplication by $2$. Hence the assumption that  $H_*(\tilde{U})$ be finitely generated implies that $2 H_*(U)$ is finitely generated as well.

Now we look at the quotient $H_*(U)/2H_*(U)$. This is just $H_*(U) \otimes \mathbb{Z}_2$, which by the Universal Coefficient theorem injects as a subgroup of $H_*(U;\mathbb{Z}_2)$. We are going to show that the latter is finitely generated. To check this we use the following exact sequence involving the maps $p_*$ and $\pi_*$ (see the section on the Borsuk-Ulam theorem in \cite{hatcher1}): \[\xymatrix{ \ldots \ar[r] & H_i(\tilde{U};\mathbb{Z}_2) \ar[r]^{p_*} & H_i(U;\mathbb{Z}_2) \ar[r] & H_{i-1}(U;\mathbb{Z}_2) \ar[r]^-{\pi_*} & \ldots }\] Since the homology of $\tilde{U}$ is finitely generated, the Universal Coefficient theorem implies that the same holds for the groups $H_i(\tilde{U};\mathbb{Z}_2)$. Now an induction starting at $i = 0$ and using repeatedly the above sequence shows that all the groups $H_i(U;\mathbb{Z}_2)$ are finitely generated. Since their elements have all order $2$ because of the choice of coefficients, the groups end up being finite. Therefore the same is true of their subgroups $H_*(U)/2H_*(U)$.

The proof concludes by putting together the information of the previous two paragraphs. Letting $H := H_*(U)$ for brevity, we have shown that $2H$ is finitely generated and $H/2H$ is finite. Let $2h_1,\ldots,2h_r$ be generators of $2H$ and $h_{r+1},\ldots,h_s$ representatives of the finitely many classes in $H/2H$. Then any $x \in H$ falls in one class of $H/2H$, so there is $j = r+1,\ldots,s$ such that $x-h_j \in 2H$, and then $x = h_j + \sum_{i=1}^r c_i 2 h_i$ for appropriate coefficients $c_i \in \mathbb{Z}$. In other words, $H$ is finitely generated by the $h$s.
\end{proof}

Equation~\eqref{thm:cohattexp} expresses the cohomology of the expelled attractors $A_{\lambda}$ in terms of the completed basin of attraction $\hat{\mathcal{A}}$, which is usually not directly observable. However, one can circumvent that difficulty as follows. Choose coefficients $R = \mathbb{Z}$ or $\mathbb{Z}_2$ according to whether $M$ (and hence also $\hat{\mathcal{A}}$) is orientable near $K$ or not. By Theorems \ref{thm:pacojaime} and \ref{thm:finitehopf} we have $H^*(\hat{\mathcal{A}};R) = \check{H}^*(K;R)$. From the preceding theorem we also know that  $H^*(\hat{\mathcal{A}} \setminus K;R) \cong \check{H}^*(A_{\lambda};R)$. The long exact sequence of the pair $(\hat{\mathcal{A}},\hat{\mathcal{A}} \setminus K)$ then reads \begin{equation} \label{eq:exseq} \xymatrix{\ldots \ar[r] & \check{H}^q(K;R) \ar[r] & \check{H}^q(A_{\lambda};R) \ar[r] & H^{q+1}(\hat{\mathcal{A}},\hat{\mathcal{A}} \setminus K;R) \ar[r] & \ldots}\end{equation} Finally, by Alexander duality $H_q(\hat{\mathcal{A}},\hat{\mathcal{A}} \setminus K;R) \cong \check{H}^{n-q}(K;R)$. Using the Universal Coefficient theorem one can then compute $H^q(\hat{\mathcal{A}},\hat{\mathcal{A}} \setminus K;R)$ in terms of this and, in favorable circumstances, use the exact sequence in Equation \eqref{eq:exseq} to extract $\check{H}^*(A_{\lambda};R)$.

The following results illustrate this procedure by computing the cohomology of the expelled attractors in some particular cases. In the classical Hopf bifurcation where $K$ is a critical point in the plane, the expelled attractors are periodic orbits and in particular are homeomorphic to $\mathbb{S}^1$. Generalized Hopf bifurcations of points were analyzed in \cite[Theorem~4.1]{MR4421505}. The following corollary extends those results to attractors $K$ that have the cohomology of a point (e.g. global attractors in Euclidean space):

\begin{corollary} \label{cor:point} If $K$ has the cohomology of a point, the expelled attractors $A_{\lambda}$ have the cohomology of the $(n-1)$--sphere with integer coefficients. 
\end{corollary}

In order to prove the corollary we need a simple auxiliary lemma:

\begin{lemma} Let $U$ be a boundariless, connected manifold and assume that $H_1(U) = 0$. Then $U$ is orientable.
\end{lemma}
\begin{proof} Suppose $U$ is not orientable and let $p : \tilde{U} \longrightarrow U$ be its orientation cover. Assuming a choice of basepoints that we omit from the notation, the subgroup $p_*(\pi_1(\tilde{U})) \subseteq \pi_1(U)$ consists of the (homotopy classes of) loops that preserve orientation. It is a normal subgroup of index $2$ and there is a well defined surjective homomorphism $\theta : \pi_1(U) \longrightarrow \pi_1(U)/p_*(\pi_1(\tilde{U})) = \mathbb{Z}_2$ which assigns the value $0$ to the orientation-preserving loops and $1$ to the orientation reversing ones. Since $\mathbb{Z}_2$ is an Abelian group, $\theta$ descends to a surjective homomorphism from the Abelianization of $\pi_1(U)$ to $\mathbb{Z}_2$. However, the Abelianization of $\pi_1(U)$ is $H_1(U)$ by Hurewicz's theorem, which by assumption is zero. Hence $\theta = 0$, which contradicts its surjectivity.
\end{proof}

\begin{proof}[Proof of Corollary \ref{cor:point}] Choose $\lambda > \lambda_0$ close enough so that $K$ is a repeller with basin of repulsion $\mathcal{R}_{\lambda}$. Since $K$ has the cohomology of a point, Theorem \ref{thm:pacojaime} implies that so does $\mathcal{R}_{\lambda}$, and by the the Universal Coefficient theorem it also has the homology of a point. In particular $\mathcal{R}_{\lambda}$ is connected and, by the Lemma above, orientable. Thus $\hat{\mathcal{A}}$ is orientable and we can work with coefficients $R = \mathbb{Z}$. Alexander duality gives $H_q(\hat{\mathcal{A}} ,\hat{\mathcal{A}} \setminus K;\mathbb{Z})  = \mathbb{Z}$ for $q = n$ and $0$ otherwise. By the Universal Coefficient theorem, $H^q(\hat{\mathcal{A}} ,\hat{\mathcal{A}} \setminus K;\mathbb{Z}) = \mathbb{Z}$ for $q = n$ and $0$ otherwise. Then the exact sequence in Equation \eqref{eq:exseq} implies that $A_{\lambda}$ has the cohomology of $\mathbb{S}^{n-1}$.
\end{proof}

Similar computations, which we omit, are also possible when $K$ has the cohomology of a $k$-sphere:

\begin{corollary}\label{cor:circle} Let the ambient dimension $n \geq 3$. If $K$ has the cohomology of $\mathbb{S}^k$ ($k \geq 1$) then the expelled attractors have the cohomology of $\mathbb{S}^k\times\mathbb{S}^{n-k-1}$. This is always valid except when $k=1$ and $M$ is not orientable along $K$, in which case it remains valid with $\mathbb{Z}_2$ coefficients.
\end{corollary}

Still in the setup of Corollary~\ref{cor:circle}, if the ambient space has dimension $n = 2$ the only interesting case is when $K$ has the cohomology of $\mathbb{S}^1$ (i.e. $k=1$). One can easily show that $K$ is a topological spine (see Section~\ref{sec:compacta_that_cannot}) of an annulus or a M\"obius band. Then \cite[Theorem~4.2]{MR4421505} applies to show that $A_{\lambda}$ have the cohomology of $\mathbb{S}^0 \times \mathbb{S}^1$ or $\mathbb{S}^1$, respectively.

From a very simplistic point of view, the results above are intuitively reasonable: since $K$ ``looks like'' a well behaved compactum, one imagines a thickening of $K$ to obtain a neighbourhood that is a manifold with boundary and the expelled compactum has the cohomology of its boundary. Now we show that in many circumstances there is actually a way of finding $\check{H}^*(A_{\lambda})$ just by looking at a ``thickening'' of $K$. 

We need a purely topological definition.

\begin{definition}
    Suppose $K \subseteq M$ is a compactum and $N$ is a compact manifold (possibly with boundary) that is a neighbourhood of $K$. We shall say that $K$ is a \emph{homological spine} of $N$ if the inclusion $K \subseteq N$ induces isomorphisms in \v{C}ech cohomology.
\end{definition}

The simplest picture to have in mind is that of a polyhedron (resp. a submanifold) $K$ in $\mathbb{R}^n$, because then $K$ is a spine of any one of its regular neighbourhoods $N$ (in the piecewise linear, resp. smooth, sense).

\begin{corollary} \label{cor:spines} Suppose that $K$ is a homological spine of some $N$ that is contained in $\hat{\mathcal{A}}$. Then $\check{H}^*(A_{\lambda};R) = H^*(\partial N;R)$ where $R$ is $\mathbb{Z}$, $\mathbb{Q}$ or $\mathbb{Z}_p$ with $p$ prime if $M$ is orientable along $K$ and $\mathbb{Z}_2$ otherwise.
\end{corollary}
\begin{proof} By Theorem \ref{thm:hopfbif} we only need to show that $\hat{\mathcal{A}}\setminus K$ has the cohomology of $\partial N$. By Theorem~\ref{thm:pacojaime} the pair $(\hat{\mathcal{A}},K)$ has a trivial cohomology. Since the same holds for $(N,K)$ by assumption, the pair $(\hat{\mathcal{A}}(K),N)$ has a trivial cohomology. By excision, so does $(\hat{\mathcal{A}}\setminus K,N\setminus K)$. It only remains to show that $(N\setminus K,\partial N)$ has a trivial cohomology and then the long exact sequence of the triple $(\hat{\mathcal{A}}\setminus K,N\setminus K,\partial N)$ will show that the inclusion $\partial N \subseteq \hat{\mathcal{A}}\setminus K$ induces isomorphisms in cohomology, as required.

Now we analyze $(N\setminus K,\partial N)$. After replacing $f$ with some suitable power we may assume that $N$ is positively invariant. Since $K$ is a homological spine of $N$, the pair $(N,K)$ has a trivial cohomology and therefore so does its homeomorphic copy $(f^n(N),K)$ for every $n \geq 0$. The long exact sequence of the triple $(N,f^n(N),K)$ then shows that $(N,f^n(N))$ also has a trivial cohomology, and again by excision so does $(N\setminus  {\rm int}\ f^n(N), \partial f^n(N))$. Recall that if $M$ is orientable along $K$, then $\hat{\mathcal{A}}(K)$ is also orientable and therefore so is $N$. Applying Lefschetz duality with $\mathbb{Z}$ coefficients if $\hat{\mathcal{A}}$ is orientable or with $\mathbb{Z}_2$ coefficients if not, the pair $(N \setminus {\rm int}\ f^n(N),\partial N)$ also has a trivial homology. By the compact supports property of homology, as $n \rightarrow \infty$ one has that $(N\setminus K,\partial N)$ has a trivial homology and therefore also a trivial cohomology as was to be shown.    
\end{proof}

The following result shows that for Hopf bifurcations of families of homeomorphisms of $\mathbb{R}^n$ or $\mathbb{S}^n$ we can get a description of the cohomology of the expelled attractors just in terms of the cohomology of $K$. We use the notation $\beta \check{H}^*(K)$ and $\tau \check{H}^*(K)$ to denote the free and torsion parts of $\check{H}^*(K)$, which are well defined since the latter is finitely generated.

\begin{corollary}\label{cor:duality}
    If the maps $f_{\lambda}$ are homeomorphisms and $M$ is $\mathbb{R}^n$ or $\mathbb{S}^n$, then
    \[
    \check{H}^q(A_\lambda)\cong\check{H}^q(K)\oplus\beta \check{H}^{n-q-1}(K) \oplus \tau \check{H}^{n-q}(K),\] where the cohomology groups are all reduced. 	
		In particular, the $q$th Betti number of $A_{\lambda}$ is the sum of the $q$th and $(n-q-1)$th Betti numbers of $K$.
\end{corollary}

\begin{proof}
We only prove the result for $M=\mathbb{S}^n$ since the other case follows from this one by extending the dynamics from $\mathbb{R}^n$ to $\mathbb{S}^n$ letting the point at infinity be fixed. Also, notice that since the dynamics is generated by homeomorphisms, we have $\hat{\mathcal{A}} = \mathcal{A} \subseteq \mathbb{S}^n$. We use reduced homology and cohomology throughout.

Set $K' := \mathbb{S}^n \setminus \mathcal{A}$. By Alexander duality $\check{H}^k(K') \cong H_{n-k-1}(\mathcal{A})$, so \[\beta \check{H}^k(K') \cong \beta H_{n-k-1}(\mathcal{A}) \cong \beta H^{n-k-1}(\mathcal{A}) \cong \beta \check{H}^{n-k-1}(K)\] and \[\tau \check{H}^k(K') \cong \tau H_{n-k-1}(\mathcal{A}) \cong \tau H^{n-k}(\mathcal{A}) \cong \tau \check{H}^{n-k}(K).\] Here we have used the Universal Coefficient Theorem to relate the homology and cohomology of $\mathcal{A}$ as well as the fact that $H^*(\mathcal{A}) \cong \check{H}^*(K)$.

By Theorem \ref{thm:hopfbif} we have \begin{equation} \label{eq:casoRn} \check{H}^q(A_{\lambda}) = H^q(\mathcal{A} \setminus K) = \beta H_q(\mathcal{A} \setminus K) \oplus \tau H_{q-1}(\mathcal{A} \setminus K).\end{equation} Also, Alexander duality again gives $H_q(\mathcal{A} \setminus K) \cong \check{H}^{n-q-1}(K \cup K') = \check{H}^{n-q-1}(K) \oplus \check{H}^{n-q-1}(K')$, where the last equality holds because $K$ and $K'$ are disjoint compacta. Using the results from the preceding paragraph we then have \[\beta H_q(\mathcal{A} \setminus K) \cong \beta \check{H}^{n-q-1}(K) \oplus \beta \check{H}^{n-q-1}(K') = \beta \check{H}^{n-q-1}(K) \oplus \beta \check{H}^q(K)\] and also \[\tau H_{q-1}(\mathcal{A} \setminus K) \cong \tau \check{H}^{n-q}(K) \oplus \tau \check{H}^{n-q}(K') = \tau \check{H}^{n-q}(K) \oplus \tau \check{H}^{q}(K).\] Plugging these into Equation \ref{eq:casoRn} gives the desired result.
\end{proof}

\begin{remark} 
Under the assumptions of Corollaries~\ref{cor:spines} or ~\ref{cor:duality}, the cohomology of the $A_{\lambda}$ is independent of the details of the dynamics that generate the bifurcation.   
\end{remark}

\section{Compacta that cannot undergo Hopf bifurcations} \label{sec:compacta_that_cannot}

Let us say that a compactum $K$ in some boundariless manifold $M$ \emph{admits} a (discrete-time) generalized Hopf bifurcation if there exists a continuous family of embeddings $f_{\lambda} : M \longrightarrow M$ such that $K$ undergoes a generalized Hopf bifurcation. This terminology conveys no new ideas and is just intended to shift the focus from the dynamics to the compactum itself. There is an analogous notion for continuous-time dynamics (semiflows).

Obviously, a necessary condition for a compactum $K$ to admit a Hopf bifurcation is that it can be realized as an attractor in the first place. We have seen in Proposition \ref{prop:hopfad} that this is also a sufficient condition in the continuous-time setting.  We also mentioned there that this is no longer true in the discrete-time setting. The following is an example:

\begin{example} \label{ex:solenoid} We refer to the left panel of Figure \ref{fig:solenoid}. Let $f$ be a homeomorphism of $\mathbb{R}^3$ that carries the large solid torus $T$ into its interior in the pattern shown in the figure. The well-known dyadic solenoid $K$ is the attractor determined by the trapping region $T$; i.e. $\bigcap_{n \geq 0} f^n(T)$. The right panel of the figure shows a drawing of $K$.

\begin{figure}[h]
    \centering
    \includegraphics[width=0.4\linewidth]{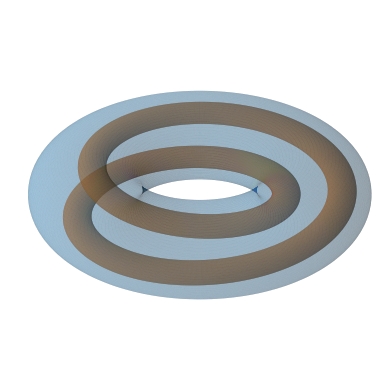}
    \includegraphics[width=0.4\linewidth]{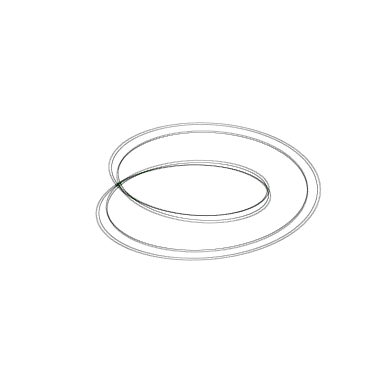}
    \caption{Construction of the dyadic solenoid as a nested sequence of solid tori.}
    \label{fig:solenoid}
\end{figure}

The continuity property of \v{C}ech cohomology shows that $\check{H}^1(K;\mathbb{Z}) \cong \mathbb{Z}[1/2]$ (the dyadic rationals). Since this group is not finitely generated, Theorem \ref{thm:finitehopf} implies that $K$ does not admit a generalized Hopf bifurcation. Thus we have an example of an attractor of a homeomorphism that does not admit a generalized Hopf bifurcation.

Of course, the same argument works for an $n$-adic solenoid ($n \geq 2$). Notice also that, since \v{C}ech cohomology is a topological invariant, no homeomorphic embedding of the solenoid in any manifold $M$ admits a Hopf bifurcation either.
\end{example}

At this point it is convenient to recall the following definition:

\begin{definition}
    A compactum $K$ in an $n$-dimensional boundariless manifold $M$ is a \emph{topological spine} if it has a compact neighbourhood $P\subseteq M$ such that there exists a homeomorphism $h : P\setminus K \longrightarrow {\rm fr}(P) \times[0,+\infty)$, where ${\rm fr}(P)$ denotes the topological frontier of $P$ in $M$. We require that $h(x) = (x,0)$ for every $x \in {\rm fr}(P)$.
\end{definition}

\begin{figure}
    \centering
    \includegraphics[width=0.5\linewidth]{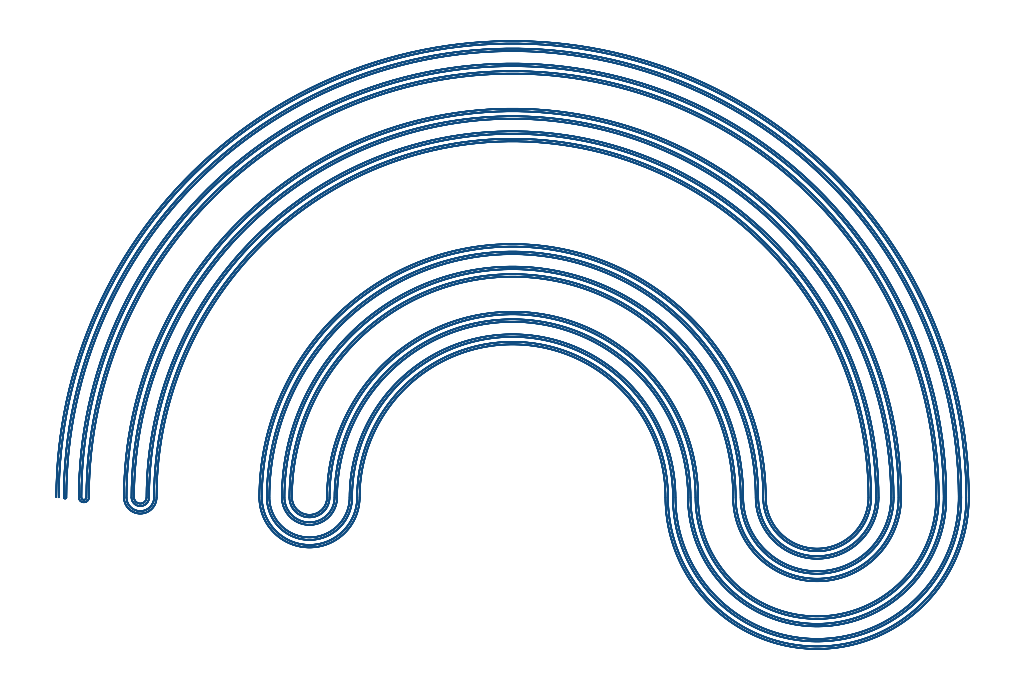}
    \caption{The Knaster continuum, a plane non-separating continuum homeomorphic to the Whitehead continuum}
    \label{fig:Knaster}
\end{figure}

\begin{figure}[htbp]
\centering
\begin{tikzpicture}[scale=1.25, >=Stealth]

  
  \fill[blue!22, opacity=0.85] (-3.2, -0.05) -- (-2.0, -0.45) -- (-2.0, -1.15) -- (-3.2, -0.75) -- cycle;
  \fill[blue!16, opacity=0.85] (-2.0, -0.45) -- (0.0, -0.20) -- (0.0, -0.90) -- (-2.0, -1.15) -- cycle;
  \fill[blue!16, opacity=0.85] (0.0, -0.20) -- (2.0, -0.45) -- (2.0, -1.15) -- (0.0, -0.90) -- cycle;
  \fill[blue!22, opacity=0.85] (2.0, -0.45) -- (3.2, -0.05) -- (3.2, -0.75) -- (2.0, -1.15) -- cycle;

  \fill[blue!40, opacity=0.75] (-1.9, 0.12) -- (-1.5, 0.28) -- (-1.1, 0.12) -- (-1.1, -0.65) -- (-1.5, -0.48) -- (-1.9, -0.65) -- cycle;
  \fill[blue!40, opacity=0.75] (1.1, 0.12) -- (1.5, 0.28) -- (1.9, 0.12) -- (1.9, -0.65) -- (1.5, -0.48) -- (1.1, -0.65) -- cycle;

  \fill[blue!8, opacity=0.85] 
    (-3.2, -0.05) -- (-2.0, 0.55) -- (0.0, 0.35) -- (2.0, 0.55) -- (3.2, -0.05) -- (2.0, -0.45) -- (0.0, -0.20) -- (-2.0, -0.45) -- cycle;

  \fill[white] (-1.9, 0.12) -- (-1.5, 0.28) -- (-1.1, 0.12) -- (-1.5, -0.02) -- cycle;
  \fill[white] (1.1, 0.12) -- (1.5, 0.28) -- (1.9, 0.12) -- (1.5, -0.02) -- cycle;

  \coordinate (C) at (0.0, -0.42);

  \draw[line width=10pt, red!20, opacity=0.6, line cap=round, line join=round] 
    (C) .. controls (0.6, -0.12) and (2.0, -0.10) .. (2.55, -0.40)
        .. controls (2.75, -0.72) and (1.5, -0.75) .. (C);

  \draw[line width=10pt, red!20, opacity=0.6, line cap=round, line join=round] 
    (C) .. controls (-0.6, -0.12) and (-2.0, -0.10) .. (-2.45, -0.40)
        .. controls (-2.65, -0.72) and (-1.5, -0.75) .. (C);

  
  \draw[ultra thick, red!85!black] 
    (C) .. controls (0.6, -0.12) and (2.0, -0.10) .. (2.55, -0.40)
        .. controls (2.75, -0.72) and (1.5, -0.75) .. (C);

  \draw[ultra thick, red!85!black] (-2.45, -0.62) -- (-2.45, -0.22);

  \draw[thick, red!85!black, samples=220, domain=0.03:0.75, declare function={f(\x) = -0.40 + 0.18*sin(360*0.12/(\x + 0.025));}] 
    plot ({-2.45 + 1.05*\x}, {f(\x)});

  \draw[ultra thick, red!85!black] 
    ({-2.45 + 1.05*0.75}, {-0.40 + 0.18*sin(360*0.12/(0.75 + 0.025))}) .. controls (-1.1, -0.12) and (-0.4, -0.22) .. (C);

  \draw[ultra thick, red!85!black] 
    (-2.45, -0.62) .. controls (-2.15, -0.72) and (-1.2, -0.73) .. (C);

  \fill[purple!95!black] (C) circle (2.2pt);

  \draw[thick, blue!80!black] 
    (-3.2, -0.05) -- (-2.0, 0.55) -- (0.0, 0.35) -- (2.0, 0.55) -- (3.2, -0.05) -- (2.0, -0.45) -- (0.0, -0.20) -- (-2.0, -0.45) -- cycle;

  \draw[thick, blue!80!black] (-1.9, 0.12) -- (-1.5, 0.28) -- (-1.1, 0.12) -- (-1.5, -0.02) -- cycle;
  \draw[thick, blue!80!black] (1.1, 0.12) -- (1.5, 0.28) -- (1.9, 0.12) -- (1.5, -0.02) -- cycle;

  \draw[thick, blue!80!black] (-1.9, 0.12) -- (-1.9, -0.65);
  \draw[thick, blue!80!black] (-1.5, 0.28) -- (-1.5, -0.48);
  \draw[thick, blue!80!black] (-1.1, 0.12) -- (-1.1, -0.65);
  \draw[blue!80!black, dashed] (-1.9, -0.65) -- (-1.5, -0.80) -- (-1.1, -0.65);

  \draw[thick, blue!80!black] (1.1, 0.12) -- (1.1, -0.65);
  \draw[thick, blue!80!black] (1.5, 0.28) -- (1.5, -0.48);
  \draw[thick, blue!80!black] (1.9, 0.12) -- (1.9, -0.65);
  \draw[blue!80!black, dashed] (1.1, -0.65) -- (1.5, -0.80) -- (1.9, -0.65);

  \draw[thick, blue!80!black] (-3.2, -0.05) -- (-3.2, -0.75);
  \draw[thick, blue!80!black] (-2.0, -0.45) -- (-2.0, -1.15);
  \draw[thick, blue!80!black] (0.0, -0.20) -- (0.0, -0.90);
  \draw[thick, blue!80!black] (2.0, -0.45) -- (2.0, -1.15);
  \draw[thick, blue!80!black] (3.2, -0.05) -- (3.2, -0.75);

  \draw[thick, blue!80!black] (-3.2, -0.75) -- (-2.0, -1.15) -- (0.0, -0.90) -- (2.0, -1.15) -- (3.2, -0.75);

  \node[blue!90!black, font=\small\bfseries] at (2.8, 0.45) {$P$};
  \node[blue!90!black, font=\small\bfseries] at (-2.8, 0.45) {$\partial P$};
  \node[red!85!black, font=\small\bfseries] at (-2.4, -0.88) {$K$};

\end{tikzpicture}
\caption{A Warsaw circle forming part of a topological spine $K$ of a handlebody $P$}
\label{fig:deeper_polish_circle_spine_3d}
\end{figure}
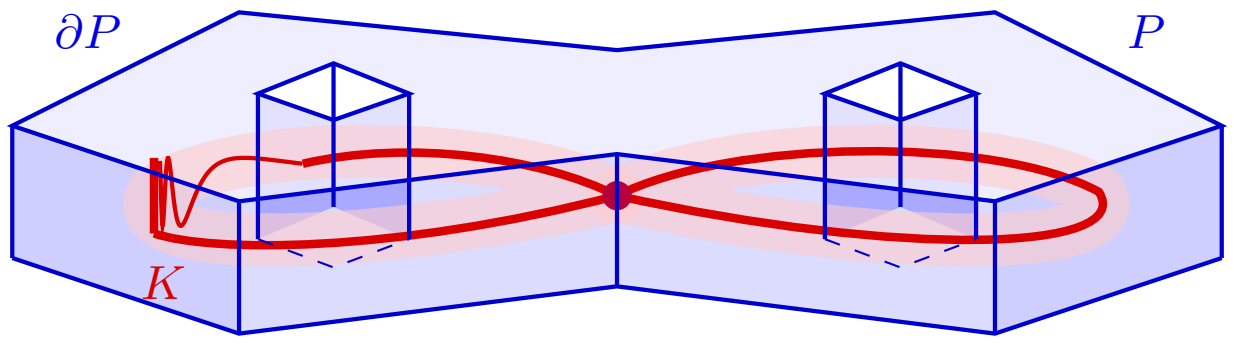

Roughly, any ``neat'' subset $K$ that can be thickened inside $M$ is a topological spine; e.g. submanifolds of smooth manifolds, subpolyhedra of PL manifolds, or subcomplexes of CW-complexes. A compactum $K \subseteq \mathbb{R}^2$ is a topological spine if and only if both $K$ and $\mathbb{R}^2 \setminus K$ have finitely many connected components. For example, the Knaster continuum in Figure \ref{fig:Knaster} is a topological spine since both $K$ and $\mathbb{R}^2 \setminus K$ are connected. Figure \ref{fig:deeper_polish_circle_spine_3d} suggests an example where the Warsaw circle is part of a spine of a genus $2$ handlebody.

The interest of topological spines in dynamics is that they are exactly those compacta that can be realized as attractors of flows. Although this is proven in \cite[Theorem~1]{gunthersegal1} or \cite[Theorem~4]{mio2} we present here a sketch of the proof for the sake of completeness.

\begin{proposition} \label{prop:spinesattractors} Let $M$ be a boundariless manifold. Any topological spine $K \subseteq M$ is an attractor of a flow in $M$. Conversely, if $K$ is an attractor of a flow (or even an injective semiflow) in $M$, it is a topological spine.
\end{proposition}
\begin{proof}[Sketch of proof] Assume $K$ is an attractor for a flow $\varphi$ and let $L$ be a Lyapunov function as in the proof of Proposition \ref{prop:hopfad}. The trajectory of any point $x \in \mathcal{A}(K) \setminus K$ intersects the set $L^{-1}(1)$ exactly once at some point $y_x = \varphi(x,t_x)$. The map $x \longmapsto t_x$ is continuous. Then $K$ is a topological spine of $P := L^{-1}([0,1])$, where a homeomorphism from $P \setminus K$ onto ${\rm fr}(P) \times [0,+\infty) = L^{-1}(1) \times [0,+\infty)$ is provided by sending a point $x \in P \setminus K$ to $(y_x,-t_x)$.

Assume now that $K$ is a topological spine of some neighbourhood $P$. There is an obvious flow in ${\rm fr}(P) \times [0,+\infty)$ that leaves every point $(x,0)$ fixed and otherwise pushes everything towards $+\infty$. This flow can be copied onto a flow in $P \setminus K$ via the homeomorphism $P \setminus K \cong {\rm fr}(P) \times [0,+\infty)$. Then it can be slowed down appropriately so it can be extended by letting every point in $K$ be fixed. The resulting flow has $K$ as an attractor.
\end{proof}

Strange attractors (e.g. the Lorenz attractor) are notorious for having a very complicated topological structure, but the proposition above tells us that even they are topological spines.

Coupling the above proposition with Proposition \ref{prop:hopfad} one gets the following purely topological characterization of what compacta admit generalized Hopf bifurcations in the continuous-time setting:

\begin{corollary} \label{cor:characterizationflows} In the context of injective semiflows, a compactum  $K$ admits a generalized Hopf bifurcation if and only if it is a topological spine.
\end{corollary}

Now we turn to the discrete-time setting. We first dispose of the case when $M$ is a surface, because then the analysis is very simple:

\begin{corollary} Let $M$ be a boundariless surface and $K \subseteq M$ a compactum. Then $K$ admits a generalized Hopf bifurcation for a family of embeddings $f_{\lambda}$ if and only if it has a finitely generated \v{C}ech cohomology.
\end{corollary}
\begin{proof} The ``only if'' part is Theorem \ref{thm:finitehopf}. The ``if'' part follows because, in a surface, a compactum with a finitely generated \v{C}ech cohomology is a topological spine of any of its frames (\cite[Theorem 6, p. 72]{moise2}).    
\end{proof}

The discussion becomes much more complicated, and one of the reasons is that now even very simple attractors may not be topological spines. Before giving examples we make the following remark:

\begin{remark} Suppose $K$ is a topological spine of some $P$ via a homeomorphism $h : P \setminus K \longrightarrow {\rm fr}(P) \times [0,+\infty)$. Then $P_n := K \cup h^{-1}({\rm fr}(P) \times [n,+\infty))$ is a nested neighbourhood basis of $K$. For $n > m$ each difference $P_m \setminus {\rm int}\ P_{n}$ is homeomorphic to the product ${\rm fr}(P) \times [m,n]$. In particular:

(1) All the $P_n$ are mutually homeomorphic.

(2) Each inclusion $P_n \subseteq P_m$ is a homotopy equivalence.

(3) The inclusion $K \subseteq P$ induces an isomorphism $\check{H}^*(P) \cong \check{H}^*(K)$ by the continuity property of \v{C}ech cohomology. Moreover, these groups are finitely generated: one can pick a compact manifold neighbourhood $N$ of $K$ contained in $P$, and since the inclusion induced isomorphism $\check{H}^*(P) \cong \check{H}^*(K)$ factors through the finitely generated cohomology of $N$, the claim follows.

\end{remark}

All these properties entail that topological spines are, actually, rather special subsets of ambient space. In fact, as mentioned above, even simple attractors for homeomorphisms may fail to be topological spines:

\begin{example} (1) The \v{C}ech cohomology of the dyadic solenoid (Example \ref{ex:solenoid}) is not finitely generated, so the solenoid is not a topological spine.

(2) Consider a modification of the construction of the solenoid, where the homeomorphism $f$ of $\mathbb{R}^3$ now carries a solid torus $T$ inside it in the pattern shown in the left panel of Figure~\ref{fig:Whitehead}. The right panel shows the continuum itself.

\begin{figure}[h]
    \centering
    \includegraphics[width=0.85\linewidth, viewport=65 170 410 355, clip]{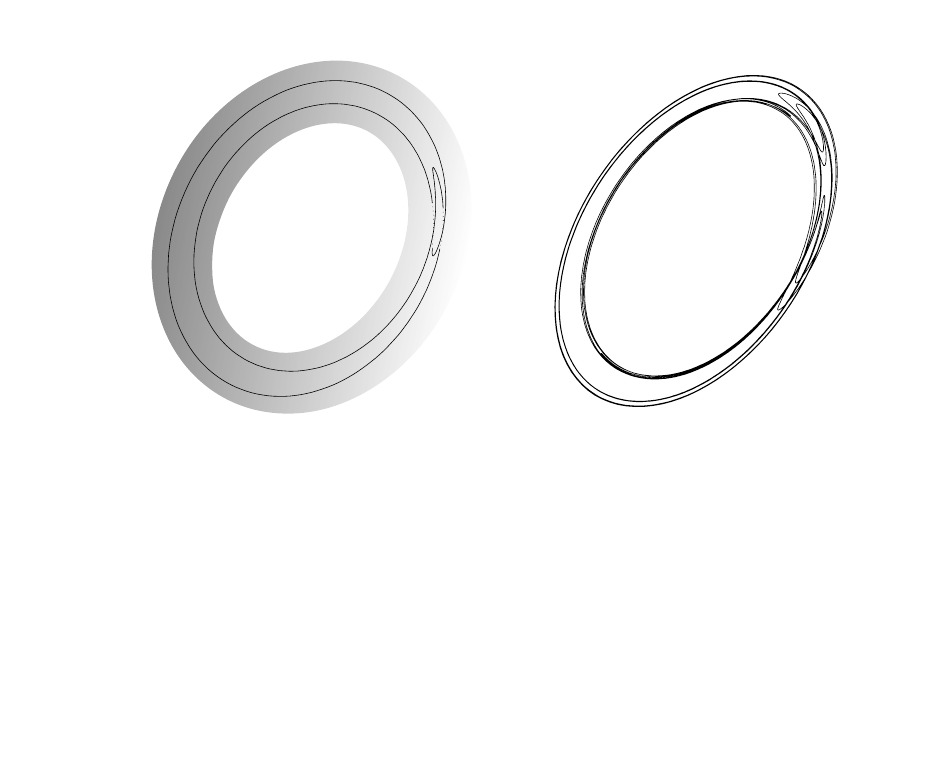}
    \caption{Construction of the Whitehead continuum with a nested sequence of solid tori}
    \label{fig:Whitehead}
\end{figure}

The Whitehead continuum is the attractor $K = \bigcap_{n \geq 0} f^n(T)$. It is not a topological spine either. For, suppose it was and call $P$ the corresponding neighbourhood. A technical argument involving generalized manifolds shows that one can take $P$ to be a tame $3$-manifold. Since each torus $f^{n+1}(T)$ is contractible inside $f^n(T)$, one easily finds that $\check{H}^*(K) = \check{H}^*({\rm point})$. Therefore, by the remarks above, $P$ would have the cohomology of a point as well. An easy argument with Lefschetz duality then shows that $\partial P$ is a $2$-sphere, and thus $P$ is a $3$-ball. Again by the remark above, this would imply that the Whitehead continuum is cellular because it has the neighbourhood basis of $3$-balls $P_n$. But it is known that the Whitehead continnuum is not cellular. Hence, it is not a topological spine. (The same argument shows that, in a $3$-dimensional manifold, any topological spine with the cohomology of a point must in fact be cellular).
\end{example}

\begin{remark} We emphasize that being a topological spine is not a topological property of the compactum: the way the compactum is embedded in ambient space is crucial. For example, the Whitehead continuum from the previous example is homeomorphic to the Knaster continuum (see \cite[Theorem~3.5, p. 40]{jubran1992chaotic}), which as mentioned earlier is a topological spine. However, the Whitehead continuum is not. This should be contrasted with the case of the dyadic solenoid: since it has a non-finitely generated \v{C}ech cohomology and the latter is a topological invariant, no homeomorphic embedding of the solenoid can ever be a topological spine. 
\end{remark}

We have just seen that two simple examples of attractors for $3$-dimensional discrete dynamical systems fail to be topological spines. We already know that the solenoid cannot undergo a generalized Hopf bifurcation. It is natural to ask if the same happens for the Whitehead continuum. Since it has the cohomology of a point, now the argument we used for the solenoid is not valid anymore. The following result establishes the somewhat surprising result that, when $M$ is a boundariless, orientable $3$-manifold, a compactum admits a Hopf bifurcation if and only if it has the rather strong property of being a topological spine:

\begin{theorem} \label{thm:onlyspines} Let $M$ be a boundariless, orientable $3$-manifold. Assume $K$ undergoes a Hopf bifurcation for a family of embeddings $f_{\lambda}$ of $M$. Then $K$ is a topological spine.
\end{theorem}

As a consequence, the Whitehead continuum does not admit a generalized Hopf bifurcation. We mention again this no-go result is a consequence of how the continuum lies in ambient space, since a homeomorphic copy of it (the Knaster continuum) does admit a bifurcation. This stands in contrast with the dyadic solenoid (or any other solenoid, in fact), for which no homeomorphic embedding admits a generalized Hopf bifurcation.

The proof of Theorem \ref{thm:onlyspines} relies heavily on results and techniques that are specific to three dimensions. Hence, we leave as an open question whether the result holds in higher dimensions:

\begin{question}\label{question}
    Let $M$ be an orientable, boundariless $n$-manifold. Is it true that if a compactum $K \subseteq M$ admits a generalized Hopf bifurcation, then $K$ is a topological spine?
\end{question}

We recall some definitions and results from $3$--manifold topology. The basic definitions about compressibility are standard and can be found for example in the classical book by Hempel \cite{hempel1}. We also assume that the objects are tame, that is, there exists a homeomorphism of the $3$--manifold that carries the object in question to some polyhedron.

We consider a $3$--manifold $B$, perhaps with boundary. By a surface $S$ in $B$ we mean a compact $2$--manifold without boundary. We always suppose that $S$ is contained in the interior of $B$. A compressing disk for $S$ (in $B$) is a disk $D \subseteq {\rm int}\ B$ that intersects $S$ precisely along $\partial D$ and such that $\partial D$ is not nullhomotopic in $S$. Observe that $2$--spheres do not have compressing disks.

Suppose $S$ is a connected surface in $B$ and $D$ is a compressing disk for $S$. Then one can compress $S$ along $D$ as follows: one removes a thin open annulus along $\partial D$ from $S$ and caps the (newly created) two boundary components of the resulting surface with two parallel copies of $D$. This produces a new closed surface $S' \subseteq {\rm int}\ B$ which satisfies $\chi(S') = \chi(S) + 2$ and may be connected or consist of two connected components. In the first case, the genus of $S'$ is one less than the genus of $S$. In the second case, both components of $S'$ have genus at least one, and their genera add up to the genus of $S$. It follows that by repeatedly compressing a surface while compressing disks are available one must arrive, in a finite number of steps, at a surface each of whose components is either a $2$-sphere or has a positive genus and no compressing disks.

For the next lemma we recall the following standard definition: a surface $S$ is said to be compressible in $B$ if it has a compressing disk $D$ or it is a $2$--sphere which bounds a $3$--ball in $B$. Otherwise, $S$ is called incompressible. 

\begin{lemma} \label{lem:goodN} Let $K$ be an attractor for an embedding $f$ of a boundariless $3$-manifold $M$. Then $K$ has a basis of neighbourhoods $N$ with the following properties:
\begin{itemize}
	\item[(i)] Each component of $\hat{\mathcal{A}}(K) \setminus K$ contains exactly one component of $\partial N$.
	\item[(ii)] Each component of $\partial N$ is incompressible in $\hat{\mathcal{A}}(K) \setminus K$.
\end{itemize}
\end{lemma}
\begin{proof} Since any $3$-manifold admits a triangulation, $K$ has a neighbourhood basis of polyhedral, compact $3$-manifolds $N$. We start with any one of them and perform the following preparatory manipulations. If necessary, we remove those components of $N$ that do not intersect $K$. Also, if a union $S_1 \cup \ldots \cup S_r$ of components of $\partial N$ bounds a compact $3$-manifold $B$ in $\hat{\mathcal{A}}(K) \setminus {\rm int}\ N$, we replace $N$ with $N \cup B$. This has the effect of removing the $D_i$ from the boundary of $\partial N$. Thus we may assume that no union of components of $\partial N$ bounds in $\hat{\mathcal{A}}(K) \setminus {\rm int}\ N$.
\smallskip

{\it Claim.} Each component of $N \setminus K$ contains exactly one component of $\partial N$.

{\it Proof of claim.} Showing that each component of $N \setminus K$ contains at least one component of $\partial N$ is an easy exercise. To prove that each component of $N \setminus K$ contains at most one component of $\partial N$ we are going to make a homological argument. Coefficients are taken in $\mathbb{Z}_2$, although we omit this from the notation. The inclusion induced map $H_2(\partial N) \longrightarrow H_2(\hat{\mathcal{A}}(K) \setminus {\rm int}\ N)$ is injective, since otherwise a union of components of $\partial N$ would bound a compact $3$-manifold in $\hat{\mathcal{A}}(K) \setminus {\rm int}\ N$, contradicting the choice of $N$. It follows from the homology exact sequence of the pair $(\hat{\mathcal{A}}(K) \setminus {\rm int}\ N, \partial N)$ that $H_3(\hat{\mathcal{A}}(K) \setminus {\rm int}\ N,\partial N) = 0$. Here we make use of $H_3(\hat{\mathcal{A}}(K) \setminus {\rm int}\ N) = 0$, which holds because every component of $\hat{\mathcal{A}}(K) \setminus{\rm int}\ N$ is a noncompact $3$-manifold. Now by excision $H_3(\hat{\mathcal{A}}(K), N) = H_3(\hat{\mathcal{A}}(K) \setminus {\rm int}\ N,\partial N) = 0$, and by the Universal Coefficient theorem, $H^3(\hat{\mathcal{A}}(K), N)=0$ as well. We know from Theorem \ref{thm:pacojaime} and Remark \ref{rem:pacojaime} that $\check{H}^2(\hat{\mathcal{A}}(K),K) = 0$, and then the cohomology long exact sequence for the triple $(\hat{\mathcal{A}}(K),N,K)$ gives $\check{H}^2(N,K) = 0$. By Alexander duality $H_1(N \setminus K,\partial N) = \check{H}^2(N,K) = 0$. It follows from the homology exact sequence for the pair $(N \setminus K,\partial N)$ that the inclusion induced map $H_0(\partial N) \longrightarrow H_0(N \setminus K)$ is injective. This means that each component of $N \setminus K$ contains at most one component of $\partial N$. $_{\blacksquare}$
\smallskip

Now we prove (i). Pick $k$ large enough so that $N \subseteq {\rm int}\ \hat{f}^{-k}(N)$. Evidently every component of $N \setminus K$ is contained in a component of $\hat{f}^{-k}(N) \setminus K$, and conversely: since every component of $\hat{f}^{-k}(N) \setminus K$ is adherent to $K$, it must intersect $N \setminus K$ and therefore contain at least one of its components. Since $\hat{f}^{-k}(N) \setminus K$ and $N \setminus K$ are homeomorphic via $\hat{f}^k$, both have the same number of components and it follows that each component of $\hat{f}^{-k}(N) \setminus K$ contains one and only one component of $N \setminus K$. Then the same is true of $\hat{\mathcal{A}}(K) \setminus K$, since it can be written as the increasing union $\bigcup_{n \geq 0} \hat{f}^{-kn}(N) \setminus K$. Coupled with the claim above, it follows that each component of $\hat{\mathcal{A}}(K) \setminus K$ contains exactly one component of $\partial N$.

Finally we show how $N$ can be modified to achieve (ii). Suppose some component $S$ of $\partial N$ has a compressing $D$ disk in $\hat{\mathcal{A}}(K) \setminus K$. Obviously $D$ is contained in the component $U$ of $\hat{\mathcal{A}}(K) \setminus K$ that contains $S$, and by (i) there is no other component of $\partial N$ in $U$, so $D \cap (\partial N) = D \cap S = \partial D$. Thus ${\rm int}\ D$ is either contained in the interior of $N$ or disjoint from $N$. In the first case we cut $N$ along $D$; in the second we attach a $2$-handle with core $D$ onto $N$. In either case we obtain a new neighbourhood $N'$ of $K$ with the same boundary as $N$ save for the component $S$, which has been compressed along $D$ to yield a new surface $S'$ which may be connected or consist of two connected components. In the latter case one of the components may bound a compact $3$-manifold $B$ in $\hat{\mathcal{A}}(K) \setminus {\rm int}\ N'$, in which case we replace $N'$ with $N' \cup B$ as we did at the beginning of the proof. By repeating this whole process finitely many times we arrive at a final stage where we get a new neighbourhood, which we call again $N$, with the property that each component of $\partial N$ has no compressing disks in $\hat{\mathcal{A}}(K) \setminus K$ and each component of $\hat{\mathcal{A}}(K) \setminus K$ still contains just one component of $\partial N$. It remains to show that every component $S$ of $\partial N$ is incompressible in $\hat{\mathcal{A}}(K) \setminus K$. If $S$ is not a $2$-sphere, this follows directly from the fact that it has no compression disks. If $S$ is a compressible $2$-sphere, it must bound a ball $B \subseteq \hat{\mathcal{A}}(K) \setminus K$. Since $B$ is entirely contained in a component of $\hat{\mathcal{A}}(K) \setminus K$ and the latter contains just one component of $\partial N$ (namely $S$), ${\rm int}\ B$ is disjoint from $\partial N$ and so it is either entirely contained in ${\rm int}\ N$ or disjoint from $N$. The latter case cannot happen since we have taken care to ensure that no union of components of $\partial N$ bounds in $\hat{\mathcal{A}}(K) \setminus {\rm int}\ N$. We now argue that we cannot have ${\rm int}\ B \subseteq {\rm int}\ N$ either. Form $DN$, the double of $N$, by pasting two copies of $N$ identifying their boundaries via the identity map. Since $\partial B = S$ is a component of $\partial N$, the double of $N$ contains the double of $B$. Both are boundariless $3$-manifolds, so the theorem on invariance of domain implies that $DB$ is an open subset of $DN$. It is obviously also a closed subset, so $DB$ is a component of $DN$. It follows that $B$ is a component of $N$. But $B$ is disjoint from $K$ by assumption, and we had discarded all components of $N$ that were disjoint from $K$. This contradiction shows that $B$ cannot be contained in $N$, and so ${\rm int}\ B$ cannot be contained in ${\rm int}\ N$. This finishes the proof that no component of $\partial N$ is compressible in $\hat{\mathcal{A}}(K) \setminus K$.
\end{proof}

The main tool we need is the Kneser--Haken finiteness theorem, that ensures that a compact manifold $B$ only admits finitely many ``essentially different'' incompressible surfaces. We use the formulation given in \cite[Theorem C, p. 494]{wright1} in a slightly simplified form:

\begin{theorem} \label{thm:kneserhaken} Let $B$ be compact orientable $3$--manifold, possibly with boundary. There exists an integer $b$ which depends only on $B$ and has the following property: any collection of pairwise disjoint, closed, connected incompressible surfaces $S_1, \ldots, S_{b+1} \subseteq {\rm int}\ B$ contains a pair of parallel surfaces.
\end{theorem}

Saying that two closed surfaces $S,S'$ are parallel in $B$ means that $S \cup S'$ is the boundary (hence also the topological frontier) of a compact $3$-manifold in ${\rm int}\ B$ that is homeomorphic to $S \times [0,1]$.

We will need the following lemma, whose proof shows why the lemma itself is not obvious:

\begin{lemma} \label{lem:paralleltransitive} Let $S_1, S_2, S_3 \subseteq {\rm int}\ B$ be incompressible, disjoint surfaces. If $S_1, S_2$ are parallel and $S_2, S_3$ are parallel, then $S_1, S_3$ are parallel.
\end{lemma}
\begin{proof} Let $P_1, P_2$ be the product manifold bounded by $S_1 \cup S_2$ and $S_2 \cup S_3$ respectively. If $S_1 \cap P_2 \neq \emptyset$ then since it is disjoint from $\partial P_2 = S_2 \cup S_3$ we have $S_1 \subseteq {\rm int}\ P_2$. Since $S_1$ is incompressible in $B$, it is also incompressible in $P_2$ and then a lemma of Haken (\cite[Lemma, p. 91]{haken1}) implies that $S_1$ is parallel to the components of $\partial P_2$, in particular to $S_3$. If $S_3 \cap P_1 \neq \emptyset$ a similar conclusion is reached: $S_3 \subseteq {\rm int}\ P_1$ and so $S_3$ is parallel to $S_1$. We are left with the case when $S_1$ is disjoint from $P_2$ and $S_3$ is disjoint from $P_1$. Then $\partial P_1 = {\rm fr}\ P_1$ is disjoint from ${\rm int}\ P_2$ and since the latter is connected, either ${\rm int}\ P_2 \subseteq P_1$ or ${\rm int}\ P_2 \cap P_1 = \emptyset$. In the first case we have $P_2 \subseteq P_1$ and so $S_3 \subseteq P_1$, contradicting our hypothesis. Thus ${\rm int}\ P_2 \cap P_1 = \emptyset$, which implies $P_2 \cap P_1 = S_2$ since $S_3$ is disjoint from $P_1$ by assumption. But then evidently $P_1 \cup P_2$ is a product manifold and $S_1, S_3$ are parallel.
\end{proof}

We can finally prove the theorem:

\begin{proof}[Proof of Theorem \ref{thm:onlyspines}] As in the first paragraph of the proof of Theorem \ref{thm:hopfbif}, we replace $M$ with $\hat{\mathcal{A}}(K)$ and the $f_{\lambda}$ by homeomorphisms $\hat{f}_{\lambda}$ of $\hat{\mathcal{A}}(K)$. This is not strictly necessary for the proof but allows for a more agile writing. It will suffice to show that $K$ is a topological spine of some $N$ in $\hat{\mathcal{A}}(K)$. Clearly, it is then also a topological spine of all the forward iterates of $N$, and since these eventually enter $\mathcal{A}(K) \subseteq M$, we will have shown that $K$ is a topological spine of some neighbourhood in $M$.

Choose $\lambda > 0$, so that $K$ is a repeller for $\hat{f}_{\lambda}$ with a basin of repulsion $\mathcal{R}_{\lambda}$. Pick a neighbourhood $N$ of $K$ as described in Lemma \ref{lem:goodN}. Acting on it with $\hat{f}$ we may assume that $N \subseteq \mathcal{R}_{\lambda}$; notice that properties (i) and (ii) of $N$ are preserved by this because $\hat{f}$ is a homeomorphism of $\hat{\mathcal{A}}(K) \setminus K$.

By Theorem \ref{thm:hopfbif} we know that $H^0(\hat{\mathcal{A}}(K) \setminus K)$ is finitely generated, which means that $\hat{\mathcal{A}}(K) \setminus K$ has finitely many connected components. Thus there exists an iterate $m_1 \geq 1$ such that $\hat{f}_{\lambda}^{m_1}$ leaves each component of $\hat{\mathcal{A}}(K) \setminus K$ invariant. Also, since $N$ is contained in $\mathcal{R}_{\lambda}$, there exists an iterate $m_2 \geq 1$ such that $N \subseteq {\rm int}\ \hat{f}_{\lambda}^{m_2}(N)$. For notational simplicity, and without any loss of generality, we replace $\hat{f}_{\lambda}$ with $\hat{f}_{\lambda}^{m_1m_2}$ and so assume directly that $\hat{f}_{\lambda}$ leaves each component of $\hat{\mathcal{A}}(K) \setminus K$ invariant and $N \subseteq {\rm int}\ \hat{f}_{\lambda}(N)$.

Let $N'$ be a compact $3$-manifold neighbourhood of $\overline{\mathcal{R}_{\lambda}}$ and define $B := N' \setminus {\rm int}\ N$. Let $b$ be the number for $B$ provided by the Kneser-Haken theorem quoted above (Theorem \ref{thm:kneserhaken}). Let $S_i$ be a component of $\partial N$. All its iterates $\hat{f}_{\lambda}(S_i), \hat{f}_{\lambda}^2(S_i), \ldots, \hat{f}^{b+1}_{\lambda}(S_i)$ lie in the interior of the compact manifold $B$. Because $\hat{f}_{\lambda}$ is a homeomorphism of $\hat{\mathcal{A}}(K) \setminus K$, property (ii) of $N$ implies that all those iterates are incompressible in $\hat{\mathcal{A}}(K) \setminus K$ and therefore also in $B$. Thus by the Kneser-Haken finiteness theorem two of them, say $\hat{f}_{\lambda}^j(S_i)$ and $\hat{f}_{\lambda}^{j'}(S_i)$ with $j' > j$, are parallel. Applying $\hat{f}_{\lambda}^{-j'}$ we see that $\hat{f}_{\lambda}^{-k(i)}(S_i)$ and $S_i$ are parallel in $\hat{\mathcal{A}}(K) \setminus K$, where we set $k(i) = j'-j$. Then applying $\hat{f}^{-k(i)}$ repeatedly we see that each consecutive pair $\hat{f}^{-mk(i)}(S_i), \hat{f}^{-(m+1)k(i)}(S_i)$ is parallel, and by Lemma \ref{lem:paralleltransitive} $S_i$ is parallel to all the iterates $\hat{f}^{-mk(i)}(S_i)$.

Set $k$ to be the product of the $k(i)$ as $S_i$ ranges over the components of $\partial N$. As just argued, each pair $S_i, f_{\lambda}^{-k}(S_i)$ is parallel in $\hat{\mathcal{A}}(K) \setminus K$. We denote by $Q_i$ the compact $3$-manifold they cobound; recall that $Q_i$ is homeomorphic to $S_i \times [0,1]$. Consider the compact $3$-manifold $C := N \setminus {\rm int}\ \hat{f}_{\lambda}^{-k}(N)$. This has boundary $\partial C = \partial N \cup \partial \hat{f}_{\lambda}^{-k}(N) = \partial N \cup \hat{f}_{\lambda}^{-k}(\partial N)$. It is disjoint from $K$ and so each component $C_i$ of $C$ is contained in a component of $\hat{\mathcal{A}}(K) \setminus K$; property (i) of $N$ then ensures that $\partial C_i$ consists of at most one component $S_i$ of $\partial N$ and at most one component of $\hat{f}^{-k}(\partial N)$. Clearly $C_i$ must have a nonempty boundary, and since both $S_i$ and $\hat{f}^{-k}(S_i)$ are contained in the same component of $\hat{\mathcal{A}}(K) \setminus K$, it follows that $\partial C_i = S_i \cup \hat{f}^{-k}(S_i)$.

We now show that $C_i = Q_i$. Notice that both are compact $3$-manifold with the same boundary, hence the same frontier. The connected set ${\rm int}\ Q_i$ is disjoint from the frontier of $C_i$, and therefore it is either entirely contained in $Q_i$ or disjoint from $C_i$. In the latter case $Q_i \cup C_i$ would be a compact $3$-manifold without boundary, and so it would have to be a whole component of the boundariless $3$-manifold $\hat{\mathcal{A}}(K)$. But every component of the latter intersects $K$, whereas $Q_i \cup C_i$ is disjoint from $K$ by construction. Therefore ${\rm int}\ Q_i \subseteq C_i$ and so $Q_i \subseteq C_i$. The same argument interchanging the roles of $Q_i$ and $C_i$ proves the other inclusion, which finally shows that $Q_i = C_i$.

The proof finishes with a routine argument. Each component $C_i$ of $C = N \setminus {\rm int}\ \hat{f}_{\lambda}^{-k}(N)$ is homeomorphic to $S_i \times [0,1]$, and so there is a homeomorphism $h : (\partial N) \times [0,1] \longrightarrow C$ carrying $(\partial N) \times 0$ onto $\partial N$ and $(\partial N) \times 1$ onto $\partial \hat{f}_{\lambda}^{-k}(N)$. Writing $N \setminus K = \bigcup_{n \geq 1} \hat{f}_{\lambda}^{-nk}(C)$ exhibits $N \setminus K$ as a union of homeomorphic copies of $(\partial N) \times [n,n+1]$ via the homeomorphism $h_n(x,t) := \hat{f}_{\lambda}^{-nk}(h(x,t-n))$, where each copy intersects the next along the part of its boundary that corresponds to $(\partial N) \times 1$. These homeomorphisms can be modified inductively as follows. We leave $h_0 (=h)$ as is and define a homeomorphism $g_0 : \partial N \longrightarrow \partial N$ by $h_1(g_0(x),0) = h_0(x,1)$. Then we replace $h_1$ with $h_1(x,t) := h_1(g_0(x),t)$. This makes $h_0 \cup h_1 : \partial N \times [0,2] \longrightarrow C \cup \hat{f}^{-k}_{\lambda}(C)$ a well defined homeomorphism. Continuing in the same fashion one ends up with a homeomorphism $(\partial N) \times [0,+\infty) \cong N \setminus K$. Thus $K$ is a topological spine of the compact $3$-manifold $N$.
\smallskip

\end{proof}

\begin{remark}
    It is not difficult to extract from the proof of Theorem~\ref{thm:onlyspines} that $\hat{\mathcal{A}}\setminus K$ is an ascending union of concentric neighbourhoods of $A_\lambda$. This implies directly Theorem~\ref{thm:cohattexp}, giving an alternative proof in dimension $3$. It is in higher dimensions, where we do not if the answer to Question~\ref{question} is affirmative, where the significance of the homological results from the previous section becomes clear.
\end{remark}

\section{Hopf bifurcations for semiflows}\label{sec:semi}

If we consider injective semiflows instead of embeddings we can obtain stronger results about topology of the expelled attractors. In this section we do this using the language of Borsuk's theory of shape. This theory has already been used to study this kind of bifurcations in some particular cases for flows \cite{MR2371787, MR3199752} and for discrete dynamical systems under some additional assumptions \cite{Barge2020BifurcationsSystems, MR4421505}. This theory is an approximate homotopy theory specially suitable to describe spaces that present some local pathologies as those that appear often as attractors of dynamical systems. Notice that spaces that have the same homotopy type are shape equivalent, and that shape theory coincides with the classical homotopy theory when considering spaces such polyhedra or manifolds.  

To give the reader a small taste of this theory to the reader we point out that if $K$ is a topological spine of some compact neighbourhood $P$ of $K$. Then, using that there exists a homeomorphism $h:P\setminus K\longrightarrow{\rm fr}(P)\setminus K$ with $h|_{{\rm fr}(P)}(x)=(x,0)$, we can consider the nested sequence of compacta $\{P_k\}_{k\geq 0}$ with $P_k=h^{-1}({\rm fr}(P)\times[k,\infty))\cup K$. It is clear from the construction that each $P_k$ is a neighbourhood of $K$, that for each $k\geq 0$, $P_{k+1}$ is a strong deformation retract of $P_k$ and that 
\[
K=\bigcap_{k\geq 0} P_k.
\]
Notice that in general, $K$ is not a deformation retraction of $P$ (see for instance figure~\ref{fig:deeper_polish_circle_spine_3d}) but we can get deformation retractions from $P$ into arbitrarily small neighbourhoods of $K$, having a kind of ``deformation retraction in the limit''. In this case, what we get is a so-called shape deformation retraction from $P$ onto $K$ and $P$ have the same shape. For more information about the theory of shape see \cite{mardesic1} and \cite{sanjurjo2011stability} for applications to dynamical systems.

 Since an attractor $K$ of an injective semiflow is a topological spines of some dynamically defined neighbourhood $P$  by Proposition~\ref{prop:spinesattractors}, the previous construction ensures that $P$ and $K$ are of the same shape. In particular, we shall make use of the following result by Kapintanski and Rodnianski:

\begin{theorem}[Kapitanski-Rodnianski]\cite[Theorem~3.6]{kapitanski1}
    Let $K$ be an attractor of a flow. Then, the inclusion $K\subseteq\mathcal{A}(K)$ is a shape equivalence.
\end{theorem}

The construction concerning the completion of a basin of attraction carries over to the setting of injective semiflows with minor modifications. Suppose one has an attractor $K$ for an injective semiflow $\varphi$ on the boundariless manifold $M$, and denote by $\mathcal{A}(K)$ its basin of attraction. Let $f$ be the time-one map of $\varphi$. This is an embedding having $K$ as an attractor, and so we can consider the completion $\hat{\mathcal{A}}(K)$ of its basin of attraction. Recall that this is obtained by taking countably many copies $A_i$ of $\mathcal{A}(K)$ and identifying $x \in A_i$ with $f(x) \in A_{i+1}$. The semiflow $\varphi$ can be extended to a semiflow on $\hat{\mathcal{A}}(K)$, which we still denote by $\varphi$, in the obvious way: $\varphi([x \in A_i],t) := [\varphi(x,t) \in A_i]$. An easy computation with the definitions shows that $\hat{f}^k(\varphi(x,t)) = \varphi(\hat{f}^k(x),t)$ for every $x \in \hat{\mathcal{A}}(K)$, $t \geq 0$, and $k \in \mathbb{Z}$. Now, $\varphi$ can be extended to a flow $\hat{\varphi}$ on $\hat{\mathcal{A}}(K)$ as follows. Given $t \in \mathbb{R}$, choose $k \in \mathbb{Z}$ such that $k \leq t$ and set $\hat{\varphi}(x,t) := \varphi(\hat{f}^k(x),t-k)$, or equivalently $\hat{f}^k(\varphi(x,t-k))$. It is easy to show that this is independent of the choice of $k$, which in turn shows that $\hat{\varphi}$ is continuous in $(x,t)$ because locally it is of the form $\varphi(\hat{f}^k(x),t-k)$ with $k$ fixed. Finally, it satisfies the flow condition $\hat\varphi(\hat\varphi(x,t_1),t_2) = \hat\varphi(x,t_1+t_2)$: choosing $k \leq t_1, t_2$ we have $2k \leq t_1+t_2$ and so \begin{align*} \hat\varphi(x,t_1+t_2) & = \varphi(\hat{f}^{2k}(x),t_1+t_2-2k) = \varphi(\varphi(\hat{f}^{2k}(x),t_1-k),t_2-k) =\\ & =  \varphi(\hat{f}^k(\varphi(\hat{f}^k(x),t_1-k)),t_2-k)  = \hat{\varphi}(\hat{\varphi}(x,t_1),t_2).\end{align*}

The following result relates the shape of the expelled attractors from generalized Hopf bifurcations in the case of injective semiflows with the shape of the level sets $L^{-1}(\varepsilon)$ of a Lyapunov function of the attractor for $\varepsilon>0$ sufficiently small.

\begin{theorem}\label{thm:hopfbifcont}
    Let $\varphi_\lambda$, $\lambda \in [0,1]$, be a continuous family of injective semiflows on a boundariless manifold $M$. Suppose that $K$ undergoes a generalized Hopf bifurcation for $\varphi_{\lambda}$. Then the expelled attractors have the shape of the level sets $L^{-1}(\varepsilon)$ of a Lyapunov function of $K$ for $\varepsilon>0$ sufficiently small.    
\end{theorem}
\begin{proof}    Consider the family of embeddings $f_\lambda=\varphi_\lambda(\cdot,1)$. As always, we denote by $K_{\lambda}$ a continuation through attractors of the initial attractor $K$. Let $\hat{\mathcal{A}}(K_{\lambda})$ be the completion of the basin of attraction of $K_{\lambda}$, endowed with the flow $\hat{\varphi}_{\lambda}$. Then $A_{\lambda}$ is an attractor with basin of attraction $\hat{\mathcal{A}}(K_{\lambda}) \setminus K$, and therefore the inclusion $A_{\lambda} \subseteq \hat{\mathcal{A}}(K_{\lambda}) \setminus K$ is a shape equivalence. By Corollary \ref{coro:globalcont} the latter space is homeomorphic to $\hat{\mathcal{A}}(K) \setminus K$. (Notice that the meaning of the letter $K$ in that Corollary is not the same as its meaning here). Let $L$ be a Lyapunov function for $K$ for the semiflow $\varphi_0$. One can easily check that for $\varepsilon > 0$ sufficiently small $L^{-1}(\varepsilon)$ is a section of the extended flow $\hat{\varphi}_0$ on $\hat{\mathcal{A}}(K) \setminus K$, and so there is a homeomorphism $\hat{\mathcal{A}}(K) \setminus K \cong L^{-1}(\varepsilon) \times \mathbb{R}$. In particular $\hat{\mathcal{A}}(K) \setminus K$ has the homotopy type of $L^{-1}(\varepsilon)$, and the result follows.
\end{proof}

We point out that the shape of the expelled attractors of generalized Hopf bifurcations of equilibria was already obtained in \cite{MR2371787} for flows and in \cite{MR4421505} for homeomorphisms. 

The next corollary characterizes the shape of the expelled attractors for a generalized Hopf bifurcation of any smooth manifold in the smooth case.

\begin{corollary}\label{coro:subman}
    If the  injective semiflows are smooth and $K$ is a compact smooth submanifold without boundary, then the expelled attractors have the shape of $\partial P$ where $P$ is a smooth closed tubular neighbourhood of $K$. 
\end{corollary}

\begin{proof}
    The result follows from Theorem~\ref{thm:hopfbifcont} since for a smooth injective semiflow we can choose a smooth Lyapunov $L$ for $K$ and by \cite[Corollary]{MR231409}, choosing $\varepsilon>0$ sufficiently small,  the level set $L^{-1}(\varepsilon)$ is homotopically equivalent to the boundary of a smooth closed tubular neighbourhood of $K$  and these are unique up to smooth ambient isotopy. 
\end{proof}

A direct consequence of Corollary~\ref{coro:subman} is the following generalization of \cite[Theorem~3]{MR2371787}. 

\begin{corollary}
    If $M$ has dimension $n\geq 2$, the injective semiflows are smooth and $K$ is a periodic orbit, then the expelled attractors have the shape of $\mathbb{S}^1\times\mathbb{S}^{n-2}$ if $K$ preserves orientation or $\mathbb{S}^{1}\tilde{\times}\mathbb{S}^{n-2}$ otherwise. 
\end{corollary}

\section*{Declaration of Generative AI in Scientific Writing}
The authors have used Gemini to generate TikZ code for some figures and to complete the proof of Lemma \ref{lem:recubridor}. The authors take full responsibility for the content of the publication.

\bibliographystyle{plain}
\bibliography{biblio}

\end{document}